\documentclass[11pt]{article}

  \usepackage[letterpaper, left=1in, right=1in, top=1in, bottom=1in]{geometry}
  \providecommand{\keywords}[1]{\textbf{Keywords: } #1}
  
  \usepackage[affil-it]{authblk}  % for associating author affiliations
  
  \usepackage{mdframed}
\usepackage{booktabs}
\usepackage{titlesec}
  \usepackage{enumerate}
  \usepackage{palatino}  % font face
  \usepackage{setspace}  % line spacing
  \usepackage{color}  % font color
  \usepackage[dvipsnames]{xcolor}  % adds more color options
	\usepackage[labelfont=bf]{caption}	% caption style

  \usepackage{amsmath}  % essential math package
  \usepackage{amssymb}  % math symbols
  \usepackage{amsthm}  % creating proof and theorem environments
  \usepackage{bm}  % boldface letters
  \usepackage{bbm}  % supports \mathbb-style digits, e.g. \mathbbm{1} indicating an all-one vector
  \allowdisplaybreaks[4]  % breaking multiline equations
	\usepackage{relsize}  % can increase the size of math symbols
	\usepackage{algorithm}  % floating environment for algorithmic
	\usepackage{algpseudocode}  % for easy writing of pseudo code

  \DeclareMathOperator*{\minimize}{minimize}
  \DeclareMathOperator*{\maximize}{maximize}
  \DeclareMathOperator*{\argmin}{arg\,min}
  \newcommand{\st}{\mathrm{subject\;to}}

	\theoremstyle{plain}

	\theoremstyle{definition}
	
	\newtheorem{example}{Example}

  \usepackage{pifont}  % dingbats and symbols

	\usepackage{longtable}  % for tables that span more than one page
	\usepackage{multirow}  % for columns spanning multiple rows
	\usepackage{array}  % for column specifications
	\usepackage{booktabs}  % providing nicer tables, using \toprule, \midrule, etc.
	\usepackage{lscape} % To output the table in landscape
  \usepackage{graphicx}  % enhanced support for graphics, allows e.g. \in­clude­graph­ics
  \usepackage{float}  % better positioning of figures, e.g. by using [h!]
  \usepackage{subfig}  % creating subfigures with \subfloat
\usepackage{pgf,tikz,pgfplots}
\pgfplotsset{compat=1.18}
\usepgfplotslibrary{statistics, groupplots}
\usepackage{mathrsfs}
\usepackage{epsfig}
\usepackage{epstopdf}
\usepackage{subcaption}
\usepackage{xcolor}
\usetikzlibrary{calc}

	\usepackage{xr-hyper}
	\usepackage{hyperref}
	\hypersetup{
		pdfstartview = {FitV},  % fits the width of the page to the window
		colorlinks = true,    % false: boxed links; true: colored links
		linkcolor = NavyBlue,  % color of internal links
		citecolor = NavyBlue,  % color of links to bibliography
		filecolor = NavyBlue,  % color of file links
		urlcolor = NavyBlue  % color of external links
	}  
	\makeatletter
    \newcommand*{\addFileDependency}[1]{% argument=file name and extension
    \typeout{(#1)}
    \@addtofilelist{#1}
    \IfFileExists{#1}{}{\typeout{No file #1.}}
    }
    \makeatother

    \newcommand*{\myexternaldocument}[1]{%
    \externaldocument{#1}%
    \addFileDependency{#1.tex}%
    \addFileDependency{#1.aux}%
    }
    \myexternaldocument{Final.tex}

  \usepackage[authoryear, round]{natbib}  % for natbib, use bibtex as backend
\titlespacing*{\section}{0pt}{0.8ex}{0.4ex}
\titlespacing*{\subsection}{0pt}{0.6ex}{0.3ex}
\titlespacing*{\subsubsection}{0pt}{0.5ex}{0.2ex}
\titlespacing*{\paragraph}{0pt}{0.5ex}{0.2ex}



\title{\vspace*{-2cm}An inverse mixed-integer optimization framework for learning interpretable models of expert decision making}
\author[1]{Anurag Holani}
\author[1]{Rishabh Gupta}
\author[2]{John Wassick}
\author[1]{Qi Zhang \thanks{Corresponding author (qizh@umn.edu).}}
\affil[1]{Department of Chemical Engineering and Materials Science, University of Minnesota, Minneapolis, MN 55455, USA}
\affil[2]{Department of Chemical Engineering, Carnegie Mellon University, Pittsburgh, PA 15213, USA}
\date{}

\begin{document}
% \doublespacing

\maketitle
\vspace*{-1cm}

\begin{abstract}
Understanding how experts make decisions and being able to transfer that knowledge is important, especially in complex engineering applications. It is highly valuable for training novices, improving the performance of human-machine systems, and potentially enabling fully autonomous systems that perform as well as human experts.
However, an expert’s decision-making strategy, developed through years of experience, is often not directly accessible, since the implicit preferences and decision rules involved can be difficult to specify explicitly.
This has motivated the use of observed decisions made by the expert to learn an interpretable model that captures the expert’s decision-making process.
% The decisions made by the expert are often hard to explain exactly, incorporating a number of implicit assumptions and preferences.
In this work, we develop an inverse optimization approach to jointly learn the decision-maker's preferences (or perceived cost) and the decision rules governing their choices.
We then demonstrate the general applicability of our approach using three case studies: shift assignment, product planning, and a real-world routing problem.
Across these case studies, modeling both perceived cost and decision rules leads to better predictions, highlighting the value of the proposed framework and its greater flexibility in capturing and replicating expert decision making.
% showcasing the model's strength in obtaining interpretable and accurate models of expert decision-making.
% Inverse optimization can be used to learn the preferences of an expert decision maker from their choices, but is limited in it's ability to capture additional decision rules which the decision maker may be imposing. 
% In this work, we have developed a novel approach to representing the decision making process, which allows us to learn these rules, 
\end{abstract}


\section{Introduction}
Decision making plays a central role in almost all human activities, ranging from everyday choices to high-stakes decisions in complex industrial environments. In contexts such as selecting consumer products, choosing medical treatments, allocating production across facilities, or routing vehicles through urban networks, decision makers must weigh multiple trade-offs involving costs, risks, timing, quality, and feasibility. Although experts often make these decisions efficiently and with apparent ease, the underlying reasoning processes are typically only partially articulated, grounded in tacit knowledge, or constrained by organizational policies and heuristics that are not explicitly stated. In practice, what is often available to analysts and modelers is not the internal rationale of the decision maker, but the \emph{final decisions} themselves. This gives rise to the following central question: \emph{Given a set of observed decisions, can we unveil the hidden preferences, constraints, or rules that governed them?}

Understanding decision-making mechanisms is valuable in many settings. In industrial operations, for example, replicating the decision strategies of experienced production planners can lead to more consistent production schedules, more reliable resource allocation, and ultimately improved economic performance. In consumer-facing applications, uncovering the implicit preferences of users can inform product design, pricing, and personalization. In safety-critical domains, such as healthcare or emergency response, capturing the logic used by expert practitioners provides transparency, trust, and actionable guidelines for automated decision-support systems. Across these applications, the ability to infer not only \emph{what} decision makers do, but \emph{why} they do it, enables more interpretable and robust decision models.

Despite its importance, reconstructing the underlying decision-making process from observed choices is difficult for several reasons. First, individuals and organizations often cannot precisely articulate the factors involved in their decisions. A consumer selecting between models of smartphones may reference intangible aspects such as ``feel” or ``style,” which in reality encode latent trade-offs over measurable properties such as weight, durability, or screen brightness. Similarly, an industrial planner may know from experience that mixing too many product families in a single facility tends to produce operational bottlenecks, but may not express this as a formal rule. Second, decision-making procedures in competitive environments are intentionally kept proprietary, meaning that internal objectives and operational rules are unavailable to external observers. Finally, decisions are frequently influenced by unmodeled behavioral regularities, domain heuristics, and soft constraints. These considerations complicate the task of inferring decision logic, even when ample decision data are available.

In the machine learning literature, the problem of learning expert decision-making behavior from observations is referred to as \emph{imitation learning} \citep{hussein_imitation_2017, zare_survey_2024}. A common imitation learning approach is \emph{behavior cloning}, which tries to find a direct mapping of environmental states to expert actions in a supervised learning fashion. Although conceptually straightforward, behavior cloning tends to overfit as the predicted decisions may not be consistent with the expert's underlying objectives \citep{codevilla2019exploring}. \emph{Inverse reinforcement learning} (IRL) attempts to address this problem by indirectly predicting the decisions through an inferred reward function that captures the expert's preferences \citep{abbeel_apprenticeship_2004}. 
% By identifying the reward parameters for which the expert’s actions are (near-)optimal, IRL obtains a compact representation of the expert's \emph{preferences}. When deployed in new environments, IRL-based agents select actions by optimizing the learned reward, enabling stronger generalization and more principled behavior than direct imitation.
IRL commonly relies on highly expressive parametric models, such as neural networks or random forests, to represent reward functions and policy mappings; these ``black-box'' models often contain many parameters and are difficult to interpret \citep{rudin_stop_2019, guidotti_survey_2019}. As a result, they are poorly suited to domains requiring auditability, such as lending, medicine, or critical infrastructure. Even when accurate, black-box models typically fail to provide explanations that align with domain knowledge or organizational procedures. Thus, while IRL methods excel at recovering objectives (or preferences), interpretability of the objectives and the ability to recover structural elements of decision-making logic remain limited.

To address the interpretability gap, substantial work has explored \emph{logical} or \emph{symbolic} models of decision making, including decision trees, rule sets, and decision lists \citep{becker_decision_2023, de_ville_decision_2013, qiao_learning_2021}. These models represent decisions through human-readable logical conditions of the form ``if condition A holds, then take action B.'' Decision trees structure these conditions hierarchically, while rule sets represent them as unordered collections of logical rules. Decision lists impose an ordering over rules, enforcing a ``first-match semantics.” These models offer significantly greater transparency than black-box machine learning methods. However, they also come with limitations. Decision trees may become excessively deep when the feature space grows, reducing interpretability and introducing instability \citep{balcan_learning_2024}. Rule sets may include overlapping or contradictory rules, requiring conflict-resolution heuristics \citep{molnar_interpretable_2019, lakkaraju_interpretable_2016}. Decision lists impose rigid sequential structures that do not naturally capture parallel or additive influences of conditions. More importantly, these logical models are designed primarily for classification or discrete decision-making tasks; they scale poorly to the mixed discrete-continuous decision spaces common in many industrial problems. Moreover, logical models are well suited for uncovering constraints or decision logic, but they do not capture the objective-driven aspects of decision making.

This divide highlights a fundamental methodological gap: existing approaches tend to learn either \emph{preferences} or \emph{rules}, but not both. 
% IRL recovers implicit objectives but assumes fixed constraints. Decision trees and rule sets recover explicit logical conditions but fail to express trade-offs. 
Yet, in many real-world settings, expert decisions arise from a combination of both elements, governed by quantitative objectives subject to domain-specific, often context-dependent decision rules. A method capable of learning both components simultaneously would provide a more faithful, interpretable, and operationally meaningful reconstruction of decision logic.

In this work, we propose a general framework for learning \emph{optimization-based decision models} that simultaneously recovers both the unknown objective function and logical decision rules that act as soft, context-dependent constraints. We take an inverse optimization (IO) approach, where we model the underlying decision-making process as an optimization problem. Originating in the early 1990s, classical IO sought to recover the cost parameters of an optimization problem for which an observed decision is exactly optimal \citep{burton_instance_1992, burton_use_1994}. More recent IO research focuses on the \emph{data-driven} setting that considers multiple, potentially noisy observations, where the goal is to learn model parameters that minimize the prediction error \citep{chan_inverse_2023, keshavarz_imputing_2011, chan_generalized_2014, mohajerin_esfahani_data-driven_2018, bertsimas_data-driven_2015, aswani_inverse_2018}; this also applies in our case. IO has been applied successfully in domains such as medical treatment planning, logistics, energy, and finance, demonstrating its ability to yield interpretable and data-efficient models. The main innovation in our proposed approach is to represent decision rules as propositional logic statements and reformulate them as mixed-integer linear constraints so that they can be embedded in the postulated optimization model. Each rule is then associated with a learnable penalty or reward, enabling the method to identify not only \emph{which} rules are relevant but also \emph{how strongly} they influence decisions in different contexts. By incorporating these rule representations into a classical IO formulation with unknown cost parameters, we obtain a unified and inherently interpretable model capable of learning both preference structures and decision logic from observed decisions.

The remainder of this paper is organized as follows. In Section \ref{sec:forward_problem}, we introduce the forward optimization problem with the proposed general decision rule structure. In Section \ref{sec:inverse_milp}, we formulate the resulting inverse mixed-integer optimization problem and develop a cutting-plane algorithm to solve it. Section \ref{sec:case_studies} presents three computational case studies considering a shift assignment problem, a production planning problem, and the Amazon Last Mile Routing Research Challenge to assess to what extent learning decision rules alongside objectives can improve prediction accuracy and reveal interpretable structures about expert decision making. Finally, we close with some concluding remarks in Section \ref{sec:conclusions}.

\section{Forward optimization problem with decision rules}
\label{sec:forward_problem}

In IO, we assume that the observed decisions result from solving an underlying optimization problem, which is also referred to as the \emph{forward optimization problem} (FOP). Traditionally, IO is used to recover the unknown cost parameters of a postulated FOP. In the following, we propose an FOP structure that also enables the systematic incorporation of a set of candidate decision rules; IO can then be applied to determine which of these decision rules are followed by the decision maker.

\subsection{Base FOP}

We start with a base FOP that is a mixed-integer linear program (MILP) of the following form:
\begin{subequations}
\label{eqn:baseFOP}
\begin{align}
\minimize_{\bar{x} \in \mathbb{R}^n \times \mathbb{Z}^p} \quad & \bar{c}(u)^{\top}\bar{x} \label{eqn:baseFOP_obj} \\
\st \quad & \bar{A}(u) \, \bar{x} \le \bar{b}(u), \label{eqn:baseFOP_cons}
\end{align}
\end{subequations}
where $\bar{x}$ are the decision variables and $u$ are the input parameters representing contextual information that affects the feasible region of the FOP. We assume that all constraints in \eqref{eqn:baseFOP_cons} are known as they supposedly capture hard physical constraints and operational rules. We further assume that parts of the cost vector $\bar{c}$ in \eqref{eqn:baseFOP_obj} are unknown, especially those related to hidden preferences of the human decision maker. Problem \eqref{eqn:baseFOP} is standard in the IO literature \citep{wang_cutting_2009, bodur_inverse_2022} except that $\bar{c}$ are commonly assumed to be constants whereas they can generally be functions of $u$ in our formulation.

\subsection{Rule-augmented FOP}
\label{subsec:fwd_model}

Rules are commonly expressed in terms of conditions and logical statements; hence, given a set of decision rules $\mathcal{R}$, we define each rule $r \in \mathcal{R}$ through a propositional logic statement (or proposition) denoted by $\Omega_r(\bar{x},u)$, which is a function of $\bar{x}$ and $u$ and returns a Boolean output, i.e., true or false. A proposition consists of literals (true/false statements) connected by logical operators such as $\land$ (and), $\lor$ (or), and $\implies$ (implication). In the following, we provide an illustrative example of a decision rule that can be expressed through such a proposition; more examples can be found in Section S1 of the supplementary material.

\begin{example}
\label{exp:proposition}
    Consider a multistage production scheduling problem where multiple units are available for each production stage. Let the manufacturing process consist of two stages with Units A and B for Stage 1 and Units C, D, and E for Stage 2. Say the human expert performing this scheduling task applies the following decision rule: Use Unit B in Stage 1 if and only if the production rate in Unit A, denoted by $x_A$, is greater than or equal to $h$; and if Unit B is used, use Unit D or Unit E in Stage 2. With $P_B$, $P_D$, and $P_E$ denoting the Boolean variables associated with the literals corresponding to using the respective units, this decision rule can be formulated as the following proposition:
    \begin{equation}
        ((h \leq x_A) \iff P_B) \land (P_B \implies (P_D \lor P_E)). \label{eqn:example_proposition}
    \end{equation}
\end{example}

% \begin{example}
% \label{exp:proposition}
%     Consider a multistage production scheduling problem where multiple units are available for each production stage. Let the manufacturing process consist of two stages with Units A and B for Stage 1 and Units C, D, and E for Stage 2. Say the human expert performing this scheduling task applies the following decision rule: If Unit A is available ($P_{\bar{A}}$), use Unit A ($P_A$) in Stage 1, and if Unit A is used, use Unit D ($P_D$) or Unit E ($P_E$) in Stage 2. With $P_{\bar{A}}$, $P_A$, $P_D$, and $P_E$ denoting the Boolean variables associated with the corresponding literals, this decision rule can be formulated as the following proposition:
%     \begin{equation}
%         (P_{\bar{A}} \Rightarrow P_A) \land (P_A \Rightarrow (P_D \lor P_E)). \label{eqn:example_proposition}
%     \end{equation}
% \end{example}

Decision rules are often not strict, i.e., the decision maker may violate them when other potential benefits outweigh the heuristic preference. Therefore, we introduce for each rule $r$ a binary variable $z_r$ that can be used to control which rule is applied. We add the following constraint:
\begin{equation}
    z_r = 1 \implies \Omega_r(\bar{x},u) = \mathrm{true},
\end{equation}
which ensures that rule $r$ applies if $z_r=1$.

Considering decision rules of the proposed form, we can formulate the following rule-augmented FOP:
\begin{subequations}
\label{eqn:augFOP}
\begin{align}
\minimize_{\bar{x} \in \mathbb{R}^n \times \mathbb{Z}^p, \, z \in \{0,1\}^{|\mathcal{R}|}} \quad &
    \bar{c}(u)^{\top} \bar{x} - \sum_{r \in \mathcal{R}} d_r(u) \, z_r \\
\st \quad \quad \; &
    \bar{A}(u)\,\bar{x} \le \bar{b}(u) \label{eqn:base_cons} \\
&   \bigwedge_{r \in {\mathcal{R}}} \begin{bmatrix}  z_r = 1 \implies \Omega_r(\bar{x},u) = \mathrm{true} \end{bmatrix}, \label{eqn:logic_cons}
\end{align}
\end{subequations}
where $d_r$ denotes the reward for applying decision rule $r$. As such, the decision rules are incorporated as soft constraints, and the values of $d$ indicate the decision maker's tendency to use each one of them. This rule-augmented FOP serves as an interpretable approximation of the decision-making process: it captures latent objective trade-offs and heuristic behavioral rules in a unified framework. The inverse optimization task that follows seeks to estimate $(\bar{c}, d)$ so that the resulting model reproduces observed decisions as closely as possible.

\section{Data-driven inverse mixed-integer optimization approach}
\label{sec:inverse_milp}

In this section, we formulate the inverse optimization problem (IOP) associated with the rule-augmented FOP \eqref{eqn:augFOP}. We will show that the FOP can be formulated as an MILP with unknown cost coefficients, and the goal is to infer these costs from observed decisions.

\subsection{MILP formulation of the rule-augmented FOP}
\label{sec:fwd_augmented}

Recall that the decision maker’s behavior is approximated by an MILP with hard physical constraints and soft decision rules.
The decision rules, which are propositional logic statements, can be represented as mixed-integer constraints utilizing the reformulation framework of \citet{raman_modelling_1994}.
Each proposition $\Omega_r$ is a logical statement consisting of a set of literals $\mathcal{L}_r$, where we denote the Boolean variable associated with literal $\ell \in \mathcal{L}_r$ as $P_{r\ell}$. By introducing a binary variable $y_{rl}$ for each literal $P_{rl}$, where $P_{rl} = \text{true}$ if $y_{rl} = 1$ and $P_{rl} = \text{false}$ if $y_{rl} = 0$, proposition $\Omega_r$ can be represented as a set of mixed-integer linear inequalities of the following general form:
\begin{subequations}
    \begin{align}
    \hat{A}_r(u) \bar{x} + \tilde{A}_r(u) y_r \leq \hat{b}_r(u) \label{eqn:MIP_cons1} \\
    A'_r \, y_r \leq b'_r, \label{eqn:MIP_cons2}
    \end{align}
    \label{eqn:MIP_cons}
\end{subequations}
where constraints \eqref{eqn:MIP_cons1} define the relationship between $y_r$ and the variables $\bar{x}$, and constraints \eqref{eqn:MIP_cons2}, which only involve the binary variables $y_r$, are a direct reformulation of proposition $\Omega_r$.

\begin{example}
    As an example, we use the proposition from Example \ref{exp:proposition}, where we first introduce binary variables $y_{\bar{A}}$, $y_B$, $y_D$, and $y_E$ for the literals $h \leq x_A$, $P_B$, $P_D$, and $P_E$, respectively. We can then reformulate proposition \eqref{eqn:example_proposition} into the following mixed-integer linear constraints:
    \begin{equation*}
        \begin{aligned}
            h \, y_{\bar{A}} & \leq x_A \\
            x_A & \leq h + (M-h) \, y_{\bar{A}} \\
            y_{\bar{A}} & \geq y_B \\
            y_{\bar{A}} & \leq y_B \\
            y_B & \leq y_D + y_E,
        \end{aligned}
    \end{equation*}
where $M$ is a sufficiently large big-M parameter.
\end{example}

With inequalities \eqref{eqn:MIP_cons}, we can now obtain an MILP formulation of the rule-augmented FOP \eqref{eqn:augFOP} as follows:
\begin{subequations}
\label{eqn:augFOP_MILP}
\begin{align}
\minimize_{\bar{x} \in \mathbb{R}^n \times \mathbb{Z}^p, \, y \in \{0,1\}^L, \,  z \in \{0,1\}^{|\mathcal{R}|}} \quad &
    \bar{c}(u)^{\top} \bar{x} - \sum_{r \in \mathcal{R}} d_r(u) \, z_r \\
\st \quad \quad \quad \quad \; &
    \bar{A}(u)\,\bar{x} \le \bar{b}(u) \label{eqn:base_cons} \\
    & \hat{A}_r(u) \bar{x} + \tilde{A}_r(u) y_r \leq \hat{b}_r(u) \quad \forall \, r \in \mathcal{R} \\
    & A'_r \, y_r \leq b'_r + M (1 - z_r) \quad \forall \, r \in \mathcal{R}, \label{eqn:bigM}
\end{align}
\end{subequations}
where the dimensionality of $y$ is $L = \sum_{r \in \mathcal{R}} |\mathcal{L}_r|$, and $M$ is a sufficiently large big-M parameter. Note that constraints \eqref{eqn:bigM} ensure that the proposition associated with decision rule $r$ must apply if $z_r = 1$.

\subsection{IOP formulation}
\label{subsec:data_decision_loss}

For notational convenience, we write problem \eqref{eqn:augFOP_MILP} in the following compact form:
\begin{subequations}
\label{eqn:FOP}
\begin{align}
    \minimize_{x \in \mathbb{R}^{n} \times \mathbb{Z}^{m}} \quad & c(u;\theta)^{\top} x \\
    \st \quad & A(u)\,x \le b(u),
\end{align}
\end{subequations}
where we define $x := (\bar{x}, y, z)$, $c := (\bar{c}, 0, -d)$, and $A$ and $b$ such that the corresponding inequalities capture constraints \eqref{eqn:base_cons}-\eqref{eqn:bigM}. The MILP has $m = p + \sum_{r \in \mathcal{R}} |\mathcal{L}_r| + |\mathcal{R}|$ discrete variables. Moreover, we parameterize the unknown cost coefficients by $\theta$, which are to be learned from observed decisions. Here, we assume that $c(u;\theta)$ is affine in $\theta$.

We observe a finite set of contexts and decisions $\{(u_i,x_i)\}_{i \in \mathcal{I}}$, where $x_i$ denotes the decision taken under context $u_i$. A natural starting point, following data-driven IO formulations for convex FOPs \citep{aswani_inverse_2018, gupta_decomposition_2022, gupta_efficient_2023}, is to directly minimize a measure of the discrepancy between observed decisions and predicted optimal solutions. This problem can be formulated as follows:
\begin{subequations}
\label{eqn:IOP_decision}
    \begin{align}
    \minimize_{\hat{x}, \, \theta \in \Theta} \quad & \sum_{i \in \mathcal{I}} \bigl\| x_i - \hat{x}_i \bigr\| \\
    \st \quad & \hat{x}_i \in \argmin_{x \in \mathbb{R}^{n} \times \mathbb{Z}^{m}} \bigl\{ c(u_i;\theta)^{\top} x : A(u_i)\,x \le b(u_i) \bigr\} \quad \forall \, i \in \mathcal{I}, \label{eqn:lower_level}
    \end{align}
\end{subequations}
where $\hat{x}_i$ denotes the predicted decision for context $u_i$; this prediction is, according to \eqref{eqn:lower_level}, an optimal solution to the FOP under parameters $\theta$. The goal is to choose $\theta$ such that the deviation between the observed and predicted decisions is minimized. Here, we restrict $\theta$ to be selected from a given set $\Theta$.

Problem \eqref{eqn:IOP_decision} is a bilevel optimization problem with $|\mathcal{I}|$ lower-level MILPs. Due to the nonconvex nature of discrete optimization problems, one cannot directly apply a single-level KKT-based reformulation to this problem. Problem \eqref{eqn:IOP_decision} can be solved using a cutting-plane method based on a value-function reformulation \citep{wang_cutting_2009}; however, this approach introduces bilinear terms and still involves mixed-integer constraints (as detailed in Section S2 of the supplementary material) such that it quickly becomes prohibitively computationally intensive with increasing number of observations. This motivates us to consider an alternative, more amenable formulation.

\subsection{IOP formulation based on suboptimality loss}
\label{subsec:subopt_SIP}

To obtain a more tractable formulation, we adopt a closely related but slightly different IOP that minimizes the so-called suboptimality based loss \citep{moghaddass_inverse_2020, moghaddass_inverse_2021} instead of the decision loss in problem \eqref{eqn:IOP_decision}. Given parameters $\theta$, the suboptimality loss $\Delta_i(\theta) $ measures the degree to which observed decision $x_i$ is suboptimal:
\[
\Delta_i(\theta) 
:= c(u_i;\theta)^{\top} x_i - \min_{\tilde{x} \in \mathcal{S}_i} c(u_i;\theta)^{\top} \tilde{x},
\]
where
\[
\mathcal{S}_i = \{x \in \mathbb{R}^{n} \times \mathbb{Z}^{m}: A(u_i)\,x \le b(u_i)\}
\]
is the feasible set of the FOP under context $u_i$, which is assumed to be non-empty and bounded. Further, it is assumed that every observation $x_i$ is feasible, i.e., $x_i \in \mathcal{S}_i$. Then, by construction, $\Delta_i(\theta) \ge 0$ and $\Delta_i(\theta) = 0$ if and only if $x_i$ is optimal to the FOP with $\theta$ and under $u_i$.

The suboptimality-based IOP seeks to find a parameter vector $\theta$ that minimizes the total suboptimality loss. Notice that for every $\tilde{x} \in \mathcal{S}_i$, we have $c(u_i;\theta)^{\top} x_i - c(u_i;\theta)^{\top} \tilde{x} \leq \Delta_i(\theta)$. As a result, the suboptimality-based IOP can be formulated as follows:
\begin{subequations}
\label{eqn:IOP_SIP}
\begin{align}
    \minimize_{\theta \in \Theta, \, \epsilon \in \mathbb{R}^{|\mathcal{I}|}} \quad & \|\epsilon\| \\
    \st \quad &
        c(u_i;\theta)^{\top} x_i - c(u_i;\theta)^{\top} \tilde{x} \le \epsilon_i
        \quad \forall \, \tilde{x} \in \mathcal{S}_i, \, i \in \mathcal{I}, 
\end{align}
\end{subequations}
where the auxiliary variable $\epsilon_i$ takes the value of $\Delta_i(\theta)$ at the optimal solution.

\subsection{Cutting-plane algorithm}
\label{subsec:cutting_plane}

We use a cutting-plane method to solve problem \eqref{eqn:IOP_SIP}, which is generally a semi-infinite program given that $\mathcal{S}_i$ can contain infinite number of points. The cutting-plane algorithm alternates between solving a master problem with added cuts and solving a subproblem to generate new cuts if needed. With a finite subset $\bar{\mathcal{S}}_i \subseteq \mathcal{S}_i$, we can formulate the following master problem:
\begin{subequations}
\label{eqn:master_new}
\begin{align}
    \minimize_{\theta \in \Theta, \, \epsilon \in \mathbb{R}^{|\mathcal{I}|}} \quad & \|\epsilon\| \\
    \st \quad &
        c(u_i;\theta)^{\top} x_i - c(u_i;\theta)^{\top} \tilde{x} \le \epsilon_i
        \quad \forall \, \tilde{x} \in \bar{\mathcal{S}}_i, \, i \in \mathcal{I}, 
\end{align}
\end{subequations}
which is a finite-dimensional linear program.

Given a candidate parameter vector $\bar{\theta}$ obtained from solving master problem \eqref{eqn:master_new}, we check for violated constraints by solving, for each $i \in \mathcal{I}$, the following subproblem:
\begin{subequations}
\label{eqn:forward_check}
\begin{align}
    \minimize_{x \in \mathbb{R}^{n} \times \mathbb{Z}^{m}} \quad & c(u_i;\bar{\theta})^{\top} x \\
    \st \quad & A(u_i)\,x \le b(u_i).
\end{align}
\end{subequations}
Note that problem \eqref{eqn:forward_check} is simply the FOP under parameters $\bar{\theta}$ and context $u_i$.
If $c(u_i;\bar{\theta})^{\top} x_i > c(u_i;\bar{\theta})^{\top} \bar{x}_i^{*}$, with $\bar{x}^*_i$ denoting an optimal solution to problem \eqref{eqn:forward_check}, we add $\bar{x}^*_i$ to set $\bar{\mathcal{S}}_i$. If $c(u_i;\bar{\theta})^{\top} x_i \leq c(u_i;\bar{\theta})^{\top} \bar{x}_i^{*}$ for all $i \in \mathcal{I}$, we can terminate the algorithm and declare $\bar{\theta}$ the optimal solution to problem \eqref{eqn:IOP_SIP}.

The pseudocode for the proposed cutting-plane algorithm is shown in Algorithm~\ref{alg:cutting_plane_new}. The algorithm is guaranteed to converge in a finite number of iterations \citep{wang_cutting_2009} since the optimal solution to problem \eqref{eqn:forward_check} must lie at an extreme point of $\mathrm{conv}(\mathcal{S}_i)$, the convex hull of the feasible region $\mathcal{S}_i$, which is a bounded polyhedron. In the worst case, the algorithm enumerates through all extreme points, of which there are finitely many.

% Convergence is achieved when the suboptimality loss of the observation reaches zero (meaning that the observation is made optimal by the estimated cost vector), or when no new cuts are added (meaning that the sub-problem gives the exact same $x$ again, and the cost vector that minimizes the suboptimality has been found). 
% In the worst case, convergence can take an exponentially large (but finite) number of iterations, becoming impractical for larger problem sizes. In these cases, we terminate the algorithm when it plateaus in its accuracy, and the returned solution provides an approximate minimizer of \eqref{eqn:IOP_SIP}.
\begin{algorithm}[t]
\caption{Cutting-plane method for solving IOP \eqref{eqn:IOP_SIP}.}
\label{alg:cutting_plane_new}
\begin{algorithmic}[1]
    \State Initialize $\bar{\mathcal{S}}_i \gets \emptyset$ for all $i \in \mathcal{I}$.
    \Repeat
        \State Solve master problem \eqref{eqn:master_new} and obtain $(\bar{\theta},\bar{\epsilon})$.
        \State Set $\text{violated} \gets \text{false}$.
        \For{each $i \in \mathcal{I}$}
            \State Solve subproblem \eqref{eqn:forward_check} with parameters $\bar{\theta}$, obtain $\bar{x}_i^{*}$.
            \If{$c(u_i;\bar{\theta})^{\top} x_i > c(u_i;\bar{\theta})^{\top} \bar{x}_i^{*}$}
                \State Update $\bar{\mathcal{S}}_i \gets \bar{\mathcal{S}}_i \cup \{\bar{x}_i^{*}\}$.
                \State Set $\text{violated} \gets \text{true}$.
            \EndIf
        \EndFor
    \Until{$\text{violated} = \text{false}$ or a maximum number of iterations is reached}
    \State \textbf{return} $\bar{\theta}$ as the learned parameter vector.
\end{algorithmic}
\end{algorithm}

% In general mixed-integer settings, the feasible sets $\mathcal{S}_i$ may contain infinitely many distinct points due to the presence of continuous variables, and there is no universal guarantee that the cutting-plane procedure will terminate in finitely many iterations without additional assumptions (e.g., boundedness and a finite set of candidate extreme points). In practice, we monitor constraint violations and terminate when no violations larger than a prescribed tolerance $\delta$ are detected, or when a maximum iteration limit is reached. Under such stopping criteria, the returned solution provides an approximate minimizer of \eqref{eqn:IOP_SIP}. 

\section{Computational case studies}
\label{sec:case_studies}

We apply the proposed approach to three case studies: learning the preferences of a planner in shift assignment, modeling the operator decisions in production planning, and modeling the routing preferences of delivery drivers in a real-world scenario.
The case studies were implemented in Julia 1.11.2 \citep{Bezanson2017}, modeled using JuMP 1.26.0 \citep{DunningHuchetteLubin2017}, and solved with Gurobi 12.0.2 \citep{gurobi} on AMD 7763 processors of the Agate cluster of the Minnesota Supercomputing Institute.

\subsection{Shift assignment case study}
\label{shift_assignment}

For the initial case study, we consider a synthetic shift assignment problem, where the planner assigns workers to their respective shifts. This case study serves to illustrate the broad modeling capability of the proposed IO approach, by learning the objective of the planner as a function of the context under which they make the decision (such as under varying demand levels).
% $\mathcal{N}$ over a total of $\mathcal{T}$ time periods.
% , aiming to ensure sufficient staffing levels while minimizing costs.
% workers to different time slots, ensuring sufficient staffing levels while minimizing costs. There are a total of $\mathcal{N}$ workers to whom the jobs can be assigned, with a total for $\mathcal{T}$ time periods over which the assignment has to be made.

\subsubsection{Problem description}

A planner needs to assign $|\mathcal{N}|$ workers to $|\mathcal{T}|$ time periods. We know the cost of assigning each worker $n$ to time period $t$, denoted by $c_{nt}$, and the maximum number of time periods a worker can work, $T^{max}$. The input parameter $u$ represents the staffing requirement specified as the minimum number of assigments.
Moreover, it is required to assign at least one worker in each time period.
The binary decision variable $x_{nt}$ indicates whether worker $n$ is assigned to time period $t$ ($x_{nt} = 1$) or not ($x_{nt} = 0$).
% Each worker has a cost, $c_n$ associated with each shift they are assigned, which represents wage costs as well as other miscellaneous costs.
% The scheduler needs to ensure that each shift has at least one worker, and that each worker is not assigned more than a maximum number, $\mathcal{T}^{max}$, of shifts.
% Further, a minimum of $u$ shifts are to be assigned in total, based on the overall demand.

Now, the obvious approach in this scenario would be for the planner to minimize the labor costs while ensuring that the constraints outlined above are met.
However, in a real-world scenario,
% long with the objective of minimizing the labor costs while fulfilling the demand, 
the planner may wish to fulfill additional objectives, such as minimizing the inventory holding and backlog costs and respecting the scheduling preferences of workers.
These decisions may occur based on the planner's understanding of typical demand profiles and the overall demand level.
We aim to recover these objectives using the proposed IO approach.

\subsubsection{Data generation}

We generate synthetic data for the case study by generating shift assignment schedules for $|\mathcal{N}| = 5$ workers and $|\mathcal{T}| = 10$ time periods, assuming that the workers cannot work for more than $T^{max} = 8$ time periods in a day. The data is generated for $u = 10,...,40$, to get planner's preferences across differing demand scenarios, which establish the \emph{ground-truth} planner behavior.

In the ground-truth model, the planner accounts for the following three competing considerations when deciding the schedule:
\begin{itemize}
    \item the total labor cost, $\sum_{n \in \mathcal{N}} \sum_{t \in \mathcal{T}} c_n \, x_{nt}$;
    \item a desired time profile of production, represented in the objective as $\sum_{t \in \mathcal{T}}\phi(t) \sum_{n \in \mathcal{N}} x_{nt}$;
    \item and scheduling preferences of the workers, incorporated as a decision rule.
\end{itemize}
% Looking at the decision rule accounting for the worker preferences, we define:
% \[
% \Omega_{p,j} : \Bigl(\sum_{i \in P} x_{ij} \le p\Bigr),
% \]

Incorporating the above considerations, the following ground-truth model is used to obtain the training data:
\begin{subequations}
    \label{eq:data_gen}
    \begin{align}
    \minimize\limits_{x, \, y} \quad 
    & \sum\limits_{n \in \mathcal{N}} \sum\limits_{t \in \mathcal{T}} c_{nt} \, x_{nt} 
      + \sum\limits_{t \in \mathcal{T}}\phi(t) \sum\limits_{n \in \mathcal{N}} x_{nt} 
      + \gamma(u) \, z 
      \label{eq:data_gen_obj} \\
    \st \quad 
    & \sum\limits_{n \in \mathcal{N}} x_{nt} \geq 1 
      \quad \forall \, t \in \mathcal{T} 
      \label{eq:data_gen_1} \\
    & \sum\limits_{t \in \mathcal{T}} x_{nt} \leq N^{\mathrm{max}} 
      \quad \forall \, n \in \mathcal{N} 
      \label{eq:data_gen_2} \\
    & \sum\limits_{n \in \mathcal{N}} \sum\limits_{t \in \mathcal{T}} x_{nt} \geq u 
      \label{eq:data_gen_3} \\
    & z = 1 \implies x_{1t'} = 0 
      \quad \forall \, t' \in \{6, \ldots, T \} 
      \label{eq:data_gen_4} \\
    & x_{nt} \in \{0, 1\} 
      \quad \forall \, n \in \mathcal{N}, t \in \mathcal{T} 
      \label{eq:data_gen_5} \\
    & z \in \{0, 1\}
      \label{eq:data_gen_6},
    \end{align}
\end{subequations}
where $\phi(t)=80\sin(\frac{2 \pi t}{10})$ represents the temporal scheduling preferences of the planner, which typically aligns with the expected demand profile. 
Worker 1 is assumed to prefer being scheduled only for time periods before $t=5$, and the associated preference is given a weight of $\gamma = -600$ if $u \le 25$ and zero otherwise, reflecting the fact that the planner will only try to meet Worker 1's request if the demand is relatively low. 
% Finally, the cost vector, describing the labor costs per worker, is $c^\top = (250, 300, 300, 350, 350)$.
The cost of employing worker $n$ for first time period is defined as $c_{n1} = 25(1+n)$. This indicates that each successive worker, with a higher index n, has a higher cost of employment. Finally, it is assumed that since the workers prefer to work during earlier time slots, rather than the late evening ones, and hence need to be paid more for later time slots. This is included in the cost matrix, by increasing the staffing cost of each consecutive time slot by $10\%$, giving the final cost, $c_{nt} = 25(1+n)(1.1)^{t-1}$.

In the objective function \eqref{eq:data_gen_obj}, $\sum_{n \in \mathcal{N}} \sum_{t \in \mathcal{T}} c_{nt} x_{nt}$ is the total labor cost, $\sum_{t \in \mathcal{T}} \phi(t) \sum_{n \in \mathcal{N}} x_{nt}$ is used to align the production schedule to the expected demand profile, and $\gamma(u) z$ is used to account for the scheduling preferences of Worker 1.
Constraint \eqref{eq:data_gen_1} ensures that each shift has at least one worker assigned, \eqref{eq:data_gen_2} ensures that the employees are not scheduled beyond their maximum hours, \eqref{eq:data_gen_3} ensures sufficient staffing to fulfill demand, and constraint \eqref{eq:data_gen_5} enforces the decision rule by ensuring that the first worker is not scheduled from shift 6 onward.
Finally, the variables are defined in \eqref{eq:data_gen_5} and \eqref{eq:data_gen_6}.

% with costs and decision rules that are dependent on the context.

% \subsubsection{Hypothesized model without decision rules}
% \begin{subequations}
%     \begin{align}
%     \minimize\limits_{x, y, z} \quad 
%         & \sum_{n \in \mathcal{N}} \sum_{t \in \mathcal{T}} c_n x_{nt} 
%           \label{eq:shift_assignment_base_obj} \\
%     \st \quad 
%         & \sum_{n \in \mathcal{N}} x_{nt} \ge 1 
%           \qquad \forall t \in \mathcal{T}
%           \label{eq:shift_assignment_base_1} \\
%         & \sum_{t \in \mathcal{T}} x_{nt} \le N^{\mathrm{max}} 
%           \qquad \forall n \in \mathcal{N}
%           \label{eq:shift_assignment_base_2} \\
%         & \sum_{n \in \mathcal{N}} \sum_{t \in \mathcal{T}} x_{nt} \ge u
%           \label{eq:shift_assignment_base_3} \\
%         & x_{nt} \in \{0,1\} 
%           \qquad \forall n \in \mathcal{N},\, t \in \mathcal{T}
%           \label{eq:shift_assignment_base_x},
%     \end{align}
% \end{subequations}
% where the objective \eqref{eq:shift_assignment_base_obj} is used to minimize the labor costs.
% Constraint \eqref{eq:shift_assignment_base_1} ensures that each employee gets at least one shift, \eqref{eq:shift_assignment_base_2} ensures that the employees are not scheduled beyond their maximum hours, \eqref{eq:shift_assignment_base_3} ensures sufficient staffing to fulfill demand, and \eqref{eq:shift_assignment_base_x} defines the decision variables.

\subsubsection{Hypothesized model}

Now, we consider the shift assignment problem from the perspective of an external analyst, who understands the general factors being considered by the planner, but doesn't have access to the implicit objective being used by them. The labor costs are a direct observable, and hence, is known exactly. 
Hence, the analyst know that the planner has a preferred time profile of production, but do not know exactly how to represent this preference.
Further, the analyst knows tat some of the workers prefer to avoid working night shifts (shift 6 onwards), and that the planner tries to respect the worker's scheduling preferences if feasible, but neglects them is the demand is too high. However, they are not aware of which workers have these preferences, or what cut-off is established by the planner, beyond which they are unable to account for the worker's preferences in their decisions.

Based on this, the analyst hypothesizes the following:
\begin{itemize}
    \item the preferred time profile of the planner can be approximated using $\hat{\phi}(t)$, a $4^{th}$ order polynomial of time ($t$), $\hat{\phi}(t) = a_4 t^4 +a_3 t^3 + a_2 t^2 + a_1 t + a_0$, 
    \item the coefficient applied for satisfaction of decision rule can be represented using $\hat{\gamma}_n(u)$. If the decision rule applies for worker $n$, the parameter is expected to take a negative value till a certain $u$, which indicates the demand level till which the worker's preference is considered by the planner, giving $\hat{\gamma}_n(u) = {-c  \text{ if } u \leq \bar{u} \text{, else } \hat{\gamma}_n(u) = 0}$. On the other hand, if the worker has no preference, the parameter will be zero at all demand level.
\end{itemize}

Thus, we get the following hypothesized model:
\begin{subequations}
    \label{eq:shift_assignment}
    \begin{align}
    \minimize\limits_{x, y, z} \quad 
        & \sum_{n \in \mathcal{N}} \sum_{t \in \mathcal{T}} c_{nt} x_{nt} 
          + \sum_{t \in \mathcal{T}} \hat{\phi}(t) \sum_{n \in \mathcal{N}} x_{nt}
          + \sum_{n \in \mathcal{N}} \hat{\gamma}_n(u) z_n
          \label{eq:shift_assignment_obj} \\
    \st \quad 
        & \sum_{n \in \mathcal{N}} x_{nt} \ge 1 
          \qquad \forall t \in \mathcal{T}
          \label{eq:shift_assignment_1} \\
        & \sum_{t \in \mathcal{T}} x_{nt} \le N^{\mathrm{max}} 
          \qquad \forall n \in \mathcal{N}
          \label{eq:shift_assignment_2} \\
        & \sum_{n \in \mathcal{N}} \sum_{t \in \mathcal{T}} x_{nt} \ge u
          \label{eq:shift_assignment_3} \\
        &   z_n = 1 \implies
                    x_{nt'} = 0 \quad 
            \forall t' \in \{6, \dots, T\},\,n \in \mathcal{N}
            \label{eq:shift_assignment_4} \\
        & x_{nt} \in \{0,1\} 
          \qquad \forall n \in \mathcal{N},\, t \in \mathcal{T}
          \label{eq:shift_assignment_x} \\
        & z_n \in \{0,1\} 
          \qquad \forall n \in \mathcal{N}
          \label{eq:shift_assignment_y}.
    \end{align}
\end{subequations}
We make two important observations regarding \eqref{eq:shift_assignment}: the only change as compared to \eqref{eq:data_gen} is in the constraint \eqref{eq:shift_assignment_4}, which is now written for all $n \in \mathcal{N}$, and the only unknown parameters in this formulation are $\hat{\phi}(t)$ and $\hat{\gamma}_n(u)$, which can be determined using IO.
% In the objective \eqref{eq:shift_assignment_obj}, $\sum_{n \in \mathcal{N}} \sum_{t \in \mathcal{T}} c_n x_{nt}$ is used to minimize the labor costs (similar to the base formulation), $\sum_{t \in \mathcal{T}} \hat{\phi}(t) \sum_{n \in \mathcal{N}} x_{nt}$ is used to align the production schedule to the expected demand profile, and $\sum_{n \in \mathcal{N}} \hat{\gamma}_n(u) y_n$ is used to account for the employee scheduling preferences.
% , and $\sum_{n \in \mathcal{N}} \sum_{t \in \mathcal{T}} \hat{\rho}_{nt} z_{nt}$ is used to identify the time period following which the employee prefers to not be scheduled.
% Constraint \eqref{eq:shift_assignment_1} ensures that each employee gets at least one shift, \eqref{eq:shift_assignment_2} ensures that the employees are not scheduled beyond their maximum hours and \eqref{eq:shift_assignment_3} ensures sufficient staffing to fulfill demand.
Constraints \eqref{eq:shift_assignment_1}, \eqref{eq:shift_assignment_2} and \eqref{eq:shift_assignment_3} are the same as before, with the added propositional logic statement, constraint \eqref{eq:shift_assignment_4} add the new variables, $z$, which indicate whether the scheduling preferences of the employees has been respected.
\subsubsection{Results}
\begin{figure}[!htbp]
    \centering
    \begin{tikzpicture}
    \pgfplotsset{
        myaxis/.style={
            width=8cm,
            height=6cm,
            grid=both,
            grid style={gray!20},
            axis lines=left,
            xlabel style={font=\normalsize},
            ylabel style={font=\normalsize, yshift=-6pt},
            tick style={black, thick},
            tick label style={font=\normalsize},
            legend style={
                draw=black,
                fill=white,
                font=\small
            }
        }
    }

    \begin{groupplot}[
        group style={
            group name=myplots,
            group size=2 by 1,
            horizontal sep=2cm
        },
        myaxis
    ]

    \nextgroupplot[
        xmin=0, xmax=10,
        ymin=-80, ymax=120,
        xlabel={$t$},
        ylabel={$\phi(t)$},
        legend style={at={(0.05,0.05)}, anchor=south west},
        clip=false
    ]
        \refstepcounter{subfigure}
        \node[anchor=south east]
        at ([xshift=-2pt,yshift=2pt]current axis.north west)
        {(\alph{subfigure})};
        \addplot[
            only marks,
            mark=*,
            mark size=3pt,
            color=blue,
            samples at={0,1,...,10}
        ]
        {80*sin(deg(2*pi*x/10))};
        \addlegendentry{Actual}

        \addplot[
            thick,
            samples=300,
            domain=0:10,
            color=orange!85!black
        ]
        {0.029247871093749406*x^4
         + 1.5565255468750132*x^3
         - 28.12224946614593*x^2
         + 104.31603695312526*x};
        \addlegendentry{Estimated}

    \nextgroupplot[
        xmin=9, xmax=41,
        ymin=-620, ymax=20,
        xlabel={$u$},
        ylabel={$\gamma(u)$},
        legend style={at={(0.05,0.95)}, anchor=north west},
        clip=false
    ]
        \refstepcounter{subfigure}
        \node[anchor=south east]
        at ([xshift=-2pt,yshift=2pt]current axis.north west)
        {(\alph{subfigure})};
        \addplot[
            only marks,
            mark=*,
            mark size=3pt,
            color=blue,
            samples at={10,11,...,40}
        ]
        {x <= 25 ? -600 : -10};
        \addlegendentry{Actual}

        \addplot[
            only marks,
            mark=*,
            mark size=3pt,
            color=orange!85!black,
            samples at={10,11,...,40}
        ]
        {x <= 25 ? -543.6244702424478 : 10};
        \addlegendentry{Estimated}

    \end{groupplot}

    \end{tikzpicture}
    \caption{Actual vs.\ estimated values of $\phi(t)$ (a) and $\gamma(u)$ (b)}
    \label{fig:phi_gamma_side_by_side}
\end{figure}
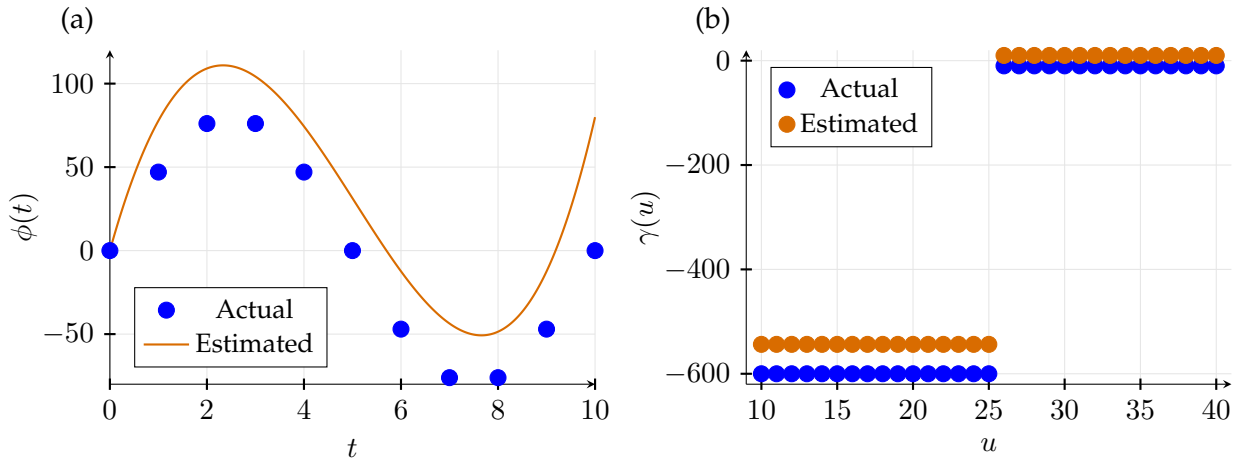

We compare the actual values of $\phi(t)$ and $\gamma_{n}(u)$ with the values estimated using IO, obtaining the trends shown in Fig. \ref{fig:phi_gamma_side_by_side} (a) and Fig. \ref{fig:phi_gamma_side_by_side} (b). We see that the predicted values of these functions correctly describes the general trends shown by the actual functions, although there is some degree of mismatch in the functional values. This is to be expected, since the regularization term in the IO formulation aims to provide the smallest parameters that can describe the behavior observed in real-world data. However, we see that the learned parameters give us a good indication of the scheduler's preferences, indicating their preferred timings of high or low production, and the degree of consideration given by them to the worker's preferences.
We also note an average normalized suboptimality loss of $0.11\%$ and an average normalized decision loss of $5.54\%$, indicating the ability of our framework to closely approximate the decisions made by the decision maker (additional details on model accuracy are provided in Section S3 of the supplementary material).
The main take-away from this case study is that our framework can determine cost parameters that are context dependent, rather than taking a fixed value under all circumstances, thus highlighting the flexibility of our IO framework.

\subsection{Production planning case study}
\label{prod_planning}

We now illustrate the proposed framework on a synthetic production planning problem.
The goal is to show that incorporating decision rules into the FOP allows us to (i) recover the decision maker's underlying objective structure
and (ii) identify rule-based preferences that strongly influence the production plan.
We also compare the prediction accuracies with and without decision rules to quantify the
impact of model misspecification.

\subsubsection{Problem description}

Consider a firm that produces a set of products $\mathcal{P}$ at a set of production facilities $\mathcal{F}$ over a planning horizon of length $H$. For each product--facility pair $(i,j) \in \mathcal{P} \times \mathcal{F}$,
we are given a production rate $r_{ij}$ (units of product $i$ per unit of time at facility $j$) and a net unit return $u_{ij}$ (dollars per unit of product $i$ produced at facility $j$).
For each product $i \in \mathcal{P}$, we have a demand $d_i$ that must be satisfied (or not exceeded) over the horizon. The decision variable $q_{ij} \ge 0$ denotes the production time allocated to
product $i$ in facility $j$. The induced quantity of product $i$ manufactured at facility $j$
is $r_{ij} \, q_{ij}$.

A standard optimization model for this problem maximizes total return across all facilities subject
to demand and capacity constraints. In practice, however, production planners may deviate
systematically from this idealized objective. For instance, they may favor products with
high unit returns, emphasize throughput, or deliberately limit the complexity of production
plans (e.g., by restricting the number of products made at each facility). Our interest is
in recovering such preferences, including rule-based complexity controls, from observed
plans.

\subsubsection{Data generation and ground-truth behavior}

We generate synthetic data for a system with $|\mathcal{P}| = 15$ products and $|\mathcal{F}| = 3$ facilities,
with a planning horizon of $H = 30$ days. The parameters $u_{ij}$ and $r_{ij}$ are set
to representative values and normalized to lie in $[0,1]$. For each instance $i \in \mathcal{I}$,
we sample a demand vector $d = (d_i)_{i \in P}$ and compute the production plan by solving
a FOP that represents the \emph{ground-truth} planner behavior.

In the ground-truth model, the planner maximizes total return, i.e., $\sum_{i \in \mathcal{P}} \sum_{j \in \mathcal{F}} u_{ij} \, r_{ij} \, q_{ij}$, and imposes a hard rule that each facility may produce at most five distinct products. Formally, the ground-truth FOP can be formulated as follows:
\begin{subequations}
\label{eqn:FOP-Prod-true}
\begin{align}
    \maximize_{q, \, x} \quad & \sum_{i \in \mathcal{P}} \sum_{j \in \mathcal{F}} u_{ij} \, r_{ij} \, q_{ij}
        \label{eqn:FOP-Prod-true-obj} \\
    \st \quad & \sum_{j \in F} r_{ij} \, q_{ij} \le d_i
        \qquad \forall \, i \in \mathcal{P}
        \label{eqn:FOP-Prod-true-demand} \\
    & \sum_{i \in P} q_{ij} \le H
        \qquad \forall \, j \in \mathcal{F}
        \label{eqn:FOP-Prod-true-cap} \\
    & h \, x_{ij} \le q_{ij} \le H \, x_{ij}
        \qquad \forall \, i \in \mathcal{P}, \, j \in \mathcal{F}
        \label{eqn:FOP-Prod-true-link} \\
    & \sum_{i \in \mathcal{P}} x_{ij} \le 5
        \qquad \forall \, j \in \mathcal{F}
        \label{eqn:FOP-Prod-true-limit} \\
    & q_{ij} \ge 0
        \qquad \forall \, i \in \mathcal{P}, \, j \in \mathcal{F}
        \label{eqn:FOP-Prod-true-qdom} \\
    & x_{ij} \in \{0,1\}
        \qquad \forall \, i \in \mathcal{P}, \, j \in \mathcal{F},
        \label{eqn:FOP-Prod-true-xdom}
\end{align}
\end{subequations}
where $x_{ij}=1$ if product $i$ is produced at facility $j$, and
$h$ is the minimum production time when a product is manufactured.
Constraints \eqref{eqn:FOP-Prod-true-limit} encode the decision rules that ``each
facility produces at most five products.'' Solving \eqref{eqn:FOP-Prod-true} for multiple demand vectors yields a collection of
observed plans $\{(d_i, (x_i, q_i))\}_{i \in \mathcal{I}}$. These observations are then used as input to the IO procedure.

\subsubsection{Hypothesized model without decision rules}

From the analyst's perspective, the ground-truth optimization problem is unknown. We assume
that the analyst knows the base constraints \eqref{eqn:FOP-Prod-true-demand}--\eqref{eqn:FOP-Prod-true-qdom}
and that the planner's objective is some linear combination of three interpretable components:

\begin{itemize}
    \item total return: $\sum_{i,j} u_{ij} r_{ij} q_{ij}$,
    \item total unit return: $\sum_{i,j} u_{ij} q_{ij}$,
    \item total production time: $\sum_{i,j} r_{ij} q_{ij}$.
\end{itemize}

Normalizing $u_{ij}$ and $r_{ij}$ to $[0,1]$, we posit the following parametric objective:
\begin{subequations}
\label{eqn:FOP-Prod}
\begin{align}
    \maximize_{q} \quad &
        \alpha \sum_{i \in P} \sum_{j \in F} u_{ij} r_{ij} q_{ij}
      + \beta  \sum_{i \in P} \sum_{j \in F} u_{ij} q_{ij}
      + \gamma \sum_{i \in P} \sum_{j \in F} r_{ij} q_{ij}
        \label{eqn:FOP-Prod-1} \\
    \st \quad &
        \sum_{j \in F} r_{ij} q_{ij} \le d_i
        \qquad \forall i \in P
        \label{eqn:FOP-Prod-2} \\
    &   \sum_{i \in P} q_{ij} \le H
        \qquad \forall j \in F
        \label{eqn:FOP-Prod-3} \\
    &   q_{ij} \ge 0
        \qquad \forall i \in P,\ \forall j \in F.
        \label{eqn:FOP-Prod-4}
\end{align}
\end{subequations}

In this \emph{no-rules} model, the parameters $(\alpha,\beta,\gamma)$ represent the relative
importance assigned by the planner to these three components. The inverse problem consists of
estimating $(\alpha,\beta,\gamma)$ from the observed plans $\{x_i\}$.

Note that \eqref{eqn:FOP-Prod} does not explicitly account for any complexity-control rule.
If the planner's behavior is strongly influenced by such a rule, the model \eqref{eqn:FOP-Prod}
is structurally misspecified, and we expect the inverse optimization procedure to struggle to
explain the data.

\subsubsection{Hypothesized model with decision rules}

To account for possible complexity-control behavior, we augment the forward model with a family
of decision rules that limit the number of products per facility. For each facility $j \in F$
and each candidate threshold $p \in \{1,\dots, |P|\}$, we define a rule:
\[
\Omega_{p,j} : \Bigl(\sum_{i \in P} x_{ij} \le p\Bigr),
\]
which is satisfied if facility $j$ produces at most $p$ distinct products. We introduce binary
variables $z_{pj} \in \{0,1\}$ indicating whether rule $\Omega_{p,j}$ is satisfied, and associate
to each rule a weight $m_{pj}$ in the objective. Negative values of $m_{pj}$ penalize violation
of the rule, consistent with the interpretation of $z_{pj}$ as yielding a reward when the rule
is satisfied.

The augmented forward problem takes the form
\begin{subequations}
\label{eqn:FOP-Prod-dec-rules}
\begin{align}
    \maximize_{q,x,y,s} \quad
        & \alpha \sum_{i \in P} \sum_{j \in F} u_{ij} r_{ij} q_{ij}
        + \beta  \sum_{i \in P} \sum_{j \in F} u_{ij} q_{ij}
        + \gamma \sum_{i \in P} \sum_{j \in F} r_{ij} q_{ij}
        + \sum_{p=1}^{|P|} \sum_{j \in F} m_{pj} z_{pj}
        \label{eqn:FOP-Prod-dec-rules-1} \\
    \st \quad &
        \sum_{j \in F} r_{ij} q_{ij} \le d_i
        \qquad \forall i \in P
        \label{eqn:FOP-Prod-dec-rules-2} \\
    &   \sum_{i \in P} q_{ij} \le H
        \qquad \forall j \in F
        \label{eqn:FOP-Prod-dec-rules-3} \\
    &   \varepsilon x_{ij} \le q_{ij} \le H x_{ij}
        \qquad \forall i \in P,\ \forall j \in F
        \label{eqn:FOP-Prod-dec-rules-4} \\
    &   \sum_{i \in P} x_{ij} = y_j
        \qquad \forall j \in F
        \label{eqn:FOP-Prod-dec-rules-5} \\
    &   z_{pj} = 1 \implies y_j \le p
        \qquad \forall p \in \{1,\dots,|P|\},\ \forall j \in F
        \label{eqn:FOP-Prod-dec-rules-6} \\
    &   q_{ij} \ge 0
        \qquad \forall i \in P,\ \forall j \in F
        \label{eqn:FOP-Prod-dec-rules-7} \\
    &   x_{ij} \in \{0,1\},\;
        y_j \in \mathbb{Z}_+,\;
        z_{pj} \in \{0,1\}
        \qquad \forall i \in P,\ \forall j \in F,\ \forall p.
        \label{eqn:FOP-Prod-dec-rules-8}
\end{align}
\end{subequations}
% where $M$ is a sufficiently large constant (e.g., $M = |P|$).
Constraint \eqref{eqn:FOP-Prod-dec-rules-6} ensures that if $z_{pj} = 1$, then $y_j \le p$.
% i.e., facility $j$ produces at most $p$ products. If $z_{pj} = 0$, the rule is effectively
% inactive. 
The cost parameters
$(\alpha,\beta,\gamma,m_{pj})$ are collected into the parameter vector $\theta$ and are
learned via inverse optimization.

In the ground-truth behavior described above, the planner enforces a maximum of five products
per facility. We therefore expect the learned rule weights $m_{5j}$ to be strongly negative
across facilities, and the weights $m_{pj}$ for other thresholds to be close to zero.

\subsubsection{Inverse optimization setup and accuracy metric}

For both the no-rules model \eqref{eqn:FOP-Prod} and the rule-augmented model
\eqref{eqn:FOP-Prod-dec-rules}, we apply the suboptimality-based inverse optimization
formulation of Section~\ref{sec:inverse_milp} to estimate the cost parameters from
synthetic data. We randomly split the generated instances into training and testing sets
and solve the inverse problem on the training set only.

To evaluate predictive performance, we compute a relative decision error on the testing set.
For each test instance $t \in \mathcal{T}$, let $x_t$ denote the observed decision vector
and $\hat{x}_t$ the optimal solution of the forward problem under the learned parameters.
We define the relative error as
\[
    e := \frac{1}{|\mathcal{T}| H_{\mathrm{tot}}}
         \sum_{t \in \mathcal{T}} \|x_t - \hat{x}_t\|_1,
\]
where $H_{\mathrm{tot}}$ is the total planning horizon aggregated across facilities
(e.g., $H_{\mathrm{tot}} = |F| H$). All results reported below are averaged over
five independent train/test splits for each training-sample size.

\subsubsection{Results with decision rules}

\begin{figure}[!htbp]
\centering

% \pgfplotscreateplotcyclelist{offon}{
%   {blue,  solid,  thick},
%   {red,   dashed, thick},
% }
\begin{tikzpicture}
\begin{axis}[
    boxplot/draw direction=y,
    ylabel={Relative error},
    yticklabel style={/pgf/number format/fixed, /pgf/number format/precision=2},
    scaled y ticks=false,
    xlabel={Number of training data points},
    height=7cm,
    boxplot={
        draw position={floor(\plotnumofactualtype/2) + 1 +(mod(\plotnumofactualtype,2)-0.5)*0.35},
        box extend=0.3,
    },
    x=1.2cm,
    boxplot/every box/.append style={solid, draw=black},
    boxplot/every whisker/.append style={solid, draw=black},
    boxplot/every median/.append style={solid, draw=black},
    boxplot/every cap/.append style={solid, draw=black},
    xtick={1,2,3,4,5,6,7},
    xticklabels={5,10,20,50,100,200,500},
    x tick label style={text width=2.5cm, align=center},
]
    \addlegendimage{area legend, draw=black, fill=cyan!100}
    \addlegendentry{Testing data predictions}
    \addlegendimage{area legend, draw=black, fill=orange!100}
    \addlegendentry{Training data predictions}

    \addplot[boxplot, draw=black, fill=orange!100]
    table[row sep=\\,y index=0] {
    data\\
    0.0022289089490694506\\
    0.07881029102452139\\
    0.13686748424674905\\
    0.046042278601412766\\
    0.0411201843538125\\
    };
    \addplot[boxplot, draw=black, fill=cyan!100]
    table[row sep=\\,y index=0] {
    data\\
    0.12223365141180642\\
    0.10593536805124253\\
    0.08106341914628805\\
    0.08036929003178424\\
    0.09154240977750895\\
    };

    \addplot[boxplot, draw=black, fill=orange!100]
    table[row sep=\\,y index=0] {
    data\\
    0.017324098485829748\\
    0.032923291073320725\\
    0.035723323183365054\\
    0.002267573696145167\\
    0.03088889934049232\\
    };
    \addplot[boxplot, draw=black, fill=cyan!100]
    table[row sep=\\,y index=0] {
    data\\
    0.05135347864767614\\
    0.07738952663971839\\
    0.05582855280415122\\
    0.049585772057555266\\
    0.07127721516513148\\
    };

    \addplot[boxplot, draw=black, fill=orange!100]
    table[row sep=\\,y index=0] {
    data\\
    0.03115578513229647\\
    0.01623205257187259\\
    0.024710540931488476\\
    0.03889680360350336\\
    0.022652929299780276\\
    };
    \addplot[boxplot, draw=black, fill=cyan!100]
    table[row sep=\\,y index=0] {
    data\\
    0.045909423793292124\\
    0.019529994701737025\\
    0.03455485643615016\\
    0.028941907023265327\\
    0.041396121511817534\\
    };

    \addplot[boxplot, draw=black, fill=orange!100]
    table[row sep=\\,y index=0] {
    data\\
    0.016006403993399587\\
    0.01814726234048763\\
    0.01106141752928001\\
    0.004684026400604392\\
    0.04411313263022873\\
    };
    \addplot[boxplot, draw=black, fill=cyan!100]
    table[row sep=\\,y index=0] {
    data\\
    0.011501840800999757\\
    0.01571897832063253\\
    0.02123322838491531\\
    0.009118365086948087\\
    0.02160073823427385\\
    };

    \addplot[boxplot, draw=black, fill=orange!100]
    table[row sep=\\,y index=0] {
    data\\
    0.00873520572962462\\
    0.01284109792050523\\
    0.009107117818345034\\
    0.014201427509960372\\
    0.0102706321785777\\
    };
    \addplot[boxplot, draw=black, fill=cyan!100]
    table[row sep=\\,y index=0] {
    data\\
    0.011499036458766067\\
    0.008541481698543492\\
    0.013174985742152763\\
    0.004401754213901284\\
    0.01083125450971616\\
    };

    \addplot[boxplot, draw=black, fill=orange!100]
    table[row sep=\\,y index=0] {
    data\\
    0.011853606558984259\\
    0.007321431132042655\\
    0.007868885728027285\\
    0.0030531527665115855\\
    0.0074503354027381065\\
    };
    \addplot[boxplot, draw=black, fill=cyan!100]
    table[row sep=\\,y index=0] {
    data\\
    0.0063886778683844465\\
    0.005925291310773274\\
    0.0070756478675325675\\
    0.005278056079647492\\
    0.005930952149669545\\
    };

    \addplot[boxplot, boxplot/draw position=6.825, draw=black, fill=orange!100]
    table[row sep=\\,y index=0] {
    data\\
    0.0032507735704344377\\
    0.009914847212046753\\
    0.005771412519678016\\
    0.006551815015052693\\
    0.006929147750310888\\
    };
    \addplot[boxplot, draw=black, fill=cyan!100]
    table[row sep=\\,y index=0] {
    data\\
    0.0032507735704344377\\
    0.0060059478002754685\\
    0.004495968062294257\\
    0.003719230553725822\\
    0.004532052023005888\\
    };

\end{axis}
\end{tikzpicture}

    \caption{Prediction accuracy on synthetic production planning data with a limit of five products per facility}
    \label{fig:5_prod_per_facility}
\end{figure}
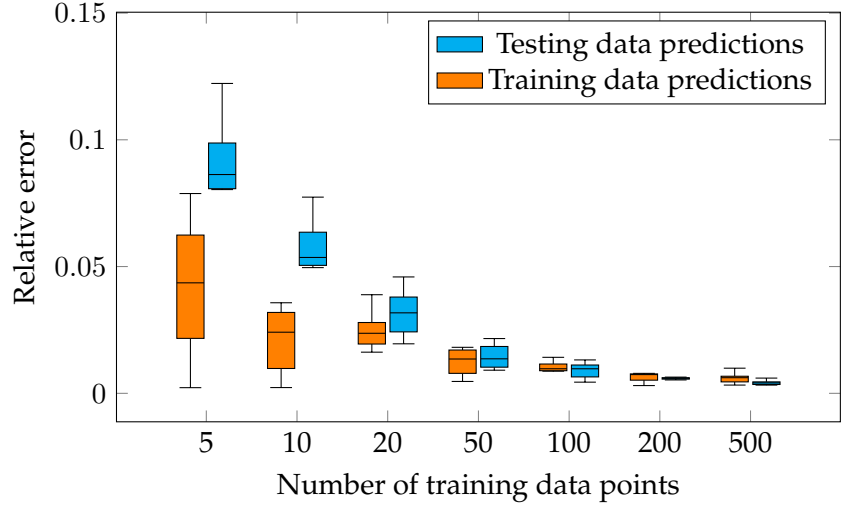

Figure~\ref{fig:5_prod_per_facility} reports the relative test error for the rule-augmented
model as a function of the number of training instances, for the case where the planner truly
limits each facility to at most five products. The prediction error decreases as more training
data are provided and stabilizes below $0.005$ for 500 training instances, indicating that the
inverse model accurately reproduces the planner's decisions.

% \begin{figure}[!htbp]
%     \centering
%     \includegraphics[scale = 0.5]{Prediction_accuracy_5_prod.png}
%     \caption{Prediction accuracy on synthetic production planning data with a limit
%     of five products per facility. The plot reports relative decision error on a held-out
%     test set as a function of the number of training instances.}
%     \label{fig:5_prod_per_facility}
% \end{figure}

The estimated objective weights for a representative run with 500 training instances are $\alpha = 48.80$, $\beta = 1.00$, and $\gamma = 1.00$.
% reported in Table~\ref{tab:obj_weights}. 
We see that $\alpha$ is substantially larger
than $\beta$ and $\gamma$, confirming that total return is the dominant objective component.

% \begin{table}[h]
% \centering
% \begin{tabular}{ccc}
% \hline
% $\alpha$ & $\beta$ & $\gamma$ \\ \hline
% 48.80 & 1.00 & 1.00 \\ \hline
% \end{tabular}
% \caption{Median estimated weights on objective components for the rule-augmented model
% (500 training instances).}
% \label{tab:obj_weights}
% \end{table}

% Table~\ref{tab:prod_limit} shows the estimated rule weights $m_{5j}$ for the three facilities.
Further, we estimated $m_{51} = -274.08$, $m_{52} = -289.46$, and $m_{51} = -205.90$.
The large negative values indicate a strong preference for satisfying the ``at most five
products per facility'' rule at each facility, consistent with the ground-truth behavior.

% \begin{table}[h]
% \centering
% \begin{tabular}{ccc}
% \hline
% $m_{51}$ & $m_{52}$ & $m_{53}$ \\ \hline
% -274.08 & -289.46 & -205.90 \\ \hline
% \end{tabular}
% \caption{Estimated rule weights for limiting the number of products to five at each facility.}
% \label{tab:prod_limit}
% \end{table}

\subsubsection{Comparison with model without decision rules}

To assess the value of decision rules, we repeat the inverse optimization procedure using
the no-rules model \eqref{eqn:FOP-Prod}. In this case, the inverse problem attempts to
explain the ground-truth plans solely through the three scalar weights $(\alpha,\beta,\gamma)$,
without any explicit representation of plan complexity.

For 500 training instances, the estimated weights were $\alpha = \beta = \gamma = 1.00$.
% the estimated weights are reported in Table~\ref{tab:obj_weights_no_dec_rules}. 
All three parameters collapsing to the same value 
indicates that the inverse problem effectively chooses a nearly flat objective to reduce
suboptimality in the presence of severe model mismatch.

% \begin{table}[!htbp]
% \centering
% \begin{tabular}{ccc}
% \hline
% $\alpha$ & $\beta$ & $\gamma$ \\ \hline
% 1.00 & 1.00 & 1.00 \\ \hline
% \end{tabular}
% \caption{Median estimated objective weights when decision rules are omitted from the
% forward model (500 training instances).}
% \label{tab:obj_weights_no_dec_rules}
% \end{table}

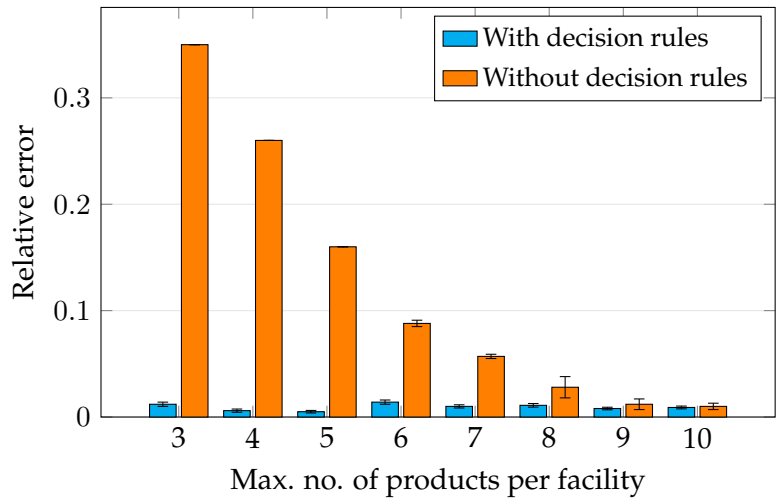
\begin{figure}[!htbp]
\centering
\begin{tikzpicture}
\begin{axis}[
    width=10.5cm,
    height=7cm,
    ybar,
    bar width=10pt,
    ymin=0,
    ylabel={Relative error},
    xlabel={Max. no. of products per facility},
    symbolic x coords={3,4,5,6,7,8,9,10},
    xtick=data,
    xtick align=inside,
    enlarge x limits=0.15,
    ymajorgrids=true,
    grid style={black!10},
    legend style={
        at={(0.98,0.98)},
        anchor=north east,
        draw=black,
        fill=white,
        font=\small,
        cells={anchor=west},
        inner xsep=3pt,
        inner ysep=2pt,
        row sep=2pt,
    },
    legend image code/.code={
      \draw[draw=black, fill]
        (0cm,-0.10cm) rectangle (0.45cm,0.10cm);
    },
    % error bars
    error bars/y dir=both,
    error bars/y explicit,
]

% ---- Series 1: With decision rules (blue) ----
\addplot[
    fill=cyan,
    draw=black,
    error bars/.cd,
    y explicit,
]
coordinates {
    (3,0.012) +- (0,0.002)
    (4,0.006) +- (0,0.0015)
    (5,0.005) +- (0,0.0012)
    (6,0.014) +- (0,0.0020)
    (7,0.010) +- (0,0.0015)
    (8,0.011) +- (0,0.0017)
    (9,0.008) +- (0,0.0012)
    (10,0.009) +- (0,0.0013)
};

% ---- Series 2: Without decision rules (orange) ----
\addplot[
    fill=orange,
    draw=black,
    error bars/.cd,
    y explicit,
]
coordinates {
    (3,0.350) +- (0,0.000)
    (4,0.260) +- (0,0.000)
    (5,0.160) +- (0,0.000)
    (6,0.088) +- (0,0.003)
    (7,0.057) +- (0,0.002)
    (8,0.028) +- (0,0.010)
    (9,0.012) +- (0,0.005)
    (10,0.010) +- (0,0.003)
};

\legend{With decision rules, Without decision rules}
\end{axis}
\end{tikzpicture}
\caption{Prediction accuracy with and without decision rules for different limits
    on the number of products that can be manufactured at a facility (500 training
    instances).}
    \label{fig:with_without_dec_rules}
\end{figure}

Figure~\ref{fig:with_without_dec_rules} compares the relative test error of the
rule-augmented model and the no-rules model for varying limits on the number of
products per facility. When the rule has a strong effect on the production plan
(e.g., low limits on the number of products), the no-rules model exhibits large
prediction errors, reflecting its inability to capture the planner's complexity
preferences. In contrast, the rule-augmented model maintains uniformly low errors
across all scenarios.

% \begin{figure}[!htbp]
%     \centering
%     \includegraphics[scale = 0.5]{Prediction_accuracy_with_without_dec_rules.png}
%     \caption{Prediction accuracy with and without decision rules for different limits
%     on the number of products that can be manufactured at a facility (500 training
%     instances).}
%     \label{fig:with_without_dec_rules}
% \end{figure}

These results demonstrate that incorporating decision rules into the forward model
is critical for accurately capturing planners' behavior when such rules materially
shape the production plans. Without decision rules, inverse optimization may produce
objective parameters that have limited behavioral interpretability and poor predictive
performance.
\subsection{Amazon last-mile routing research challenge case study}
\label{amazon_challenge}

Finally, we apply the proposed framework to a real last-mile routing problem using the
Amazon Last Mile Routing Research Challenge dataset \citep{noauthor_amazon_nodate}.
This case study serves two purposes. First, it demonstrates that our IO approach can recover interpretable decision rules from large-scale operational data.
Second, it shows that these rules substantively improve out-of-sample route prediction
relative to models that only learn pairwise travel costs.

\subsubsection{Problem description and empirical behavior}

The dataset consists of historical delivery routes executed by human drivers across
multiple depots in several U.S.\ metropolitan areas \citep{merchan_2021_2022}. Each route specifies a depot of
origin and a set of stops to be visited, along with an estimated travel-time matrix
between stops and the depot. In total, $6{,}112$ routes are provided for training and
$3{,}072$ routes for testing (aggregated across depots). The dataset also includes stop
attributes (coordinates and zone identifiers) and package attributes (e.g.,
service-time and time-window fields), although our modeling below uses only zone and
travel-time information.

A classical formulation would model each route as a Traveling Salesperson Problem (TSP)
minimizing total travel time. However, the observed driver routes deviate systematically
from shortest-time tours. Consistent with prior analyses of this dataset
\citep{zattoni_scroccaro_inverse_2024,cook_constrained_2024,guo_amazon_2023}, we observe a
strong tendency for drivers to \emph{group} stops by zone and to traverse these zones in an
organized fashion. Our objective is to learn a route-choice model that (i) captures this
zone-level structure via decision rules and (ii) refines residual preferences via IO over
pairwise zone costs.

\subsubsection{Notation and route abstraction}

Let $\mathcal{D}$ denote the set of depots and $\mathcal{R}^d$ the set of routes from depot
$d \in \mathcal{D}$. For a route $r \in \mathcal{R}^d$, let $\mathcal{S}_r$ be the set of
stops and $\mathcal{Z}_r$ the set of zones represented by stops in $\mathcal{S}_r$. We write
$\mathcal{Z}^d$ for the set of all zones observed for depot $d$ across its routes.

At the zone level, we define binary variables
\[
x_{ij} =
\begin{cases}
1, & \text{if zone } j \text{ is visited immediately after zone } i \text{ on route } r,\\
0, & \text{otherwise,}
\end{cases}
\qquad (i,j) \in \mathcal{Z}_r \times \mathcal{Z}_r,
\]
so $x$ encodes a Hamiltonian tour over $\mathcal{Z}_r$. The observed driver stop sequence
induces an observed zone tour $\tilde{x}$ by mapping each stop to its zone and compressing
contiguous runs of the same zone into a single visit. Stops missing a zone label are assigned
the zone of the nearest labeled stop on that route (in travel-time distance), ensuring every
stop belongs to exactly one zone.

\subsubsection{Baseline forward model: Restricted TSP without decision rules}

For a fixed depot $d$ and route $r \in \mathcal{R}^d$ with zone set $\mathcal{Z}_r \subseteq
\mathcal{Z}^d$, we introduce base zone-to-zone costs $\theta_{ij}^d$ for all
$(i,j) \in \mathcal{Z}^d \times \mathcal{Z}^d$. These are initialized using haversine
distances between zone centers (computed from the obfuscated coordinates) and later refined
using IO.

The baseline FOP is a TSP in the Dantzig–Fulkerson–Johnson (DFJ) form over $\mathcal{Z}_r$:
\begin{subequations}
\label{eqn:FOP-Amazon-base}
\begin{align}
    \minimize_{x} \quad &
        \sum_{i \in \mathcal{Z}^r}\sum_{j \in \mathcal{Z}^r} \theta_{ij}^d x_{ij}
        \label{eqn:FOP-Amazon-base-obj} \\
    \st \quad &
        \sum_{j \in \mathcal{Z}_r} x_{ij} = 1
        \qquad \forall i \in \mathcal{Z}_r
        \label{eqn:FOP-Amazon-base-out} \\
    &
        \sum_{i \in \mathcal{Z}_r} x_{ij} = 1
        \qquad \forall j \in \mathcal{Z}_r
        \label{eqn:FOP-Amazon-base-in} \\
    &
        \sum_{i \in \mathcal{Q}} \sum_{j \in \mathcal{Q}} x_{ij} \le |\mathcal{Q}| - 1
        \qquad \forall \emptyset \neq \mathcal{Q} \subset \mathcal{Z}_r
        \label{eqn:FOP-Amazon-base-subtour} \\
    &
        x_{ij} \in \{0,1\}
        \qquad \forall (i,j) \in \mathcal{Z}_r \times \mathcal{Z}_r
        \label{eqn:FOP-Amazon-base-binary} \\
    &
        x_{ii} = 0
        \qquad \forall i \in \mathcal{Z}_r.
        \label{eqn:FOP-Amazon-base-restrict}
\end{align}
\end{subequations}
This model captures only pairwise travel costs. As in the production-planning case, if driver
behavior is strongly shaped by higher-level decision rules, \eqref{eqn:FOP-Amazon-base} will
be structurally misspecified and IO will have limited ability to recover interpretable
preferences from sparse observations.

\subsubsection{Decision rules as zone-clustering preferences}

To represent higher-level structure, we introduce a family of decision rules that reward
transitions consistent with hypothesized zone clusters. Zone identifiers follow a structured
alphanumeric format, which we denote abstractly as $W\!-\!x.yZ$ (e.g., $P\!-\!12.3C$). We
construct candidate cluster rules by grouping zones that share selected components of these identifiers.

Let $\mathcal{K}$ be a set of candidate clustering rules. For each rule $k \in \mathcal{K}$
and zone pair $(i,j)$, we define a known indicator $M_{ij}^k \in \{0,1\}$ where $M_{ij}^k = 1$
if the transition $i \to j$ is consistent with rule $k$ (i.e., remains within the same
hypothesized cluster), and $0$ otherwise. For each depot $d$, we associate a nonnegative reward
parameter $P_k^d \ge 0$. We define the effective transition cost
\[
    c_{ij}^d(\theta^d,P^d) = \theta_{ij}^d - \sum_{k \in \mathcal{K}} P_k^d M_{ij}^k,
\]
so satisfying a rule lowers cost and encourages the route to respect the corresponding
cluster structure. Replacing $\theta^d$ by $c^d(\theta^d,P^d)$ in
\eqref{eqn:FOP-Amazon-base} yields the rule-augmented forward problem:
\begin{equation}
\label{eqn:FOP-Amazon-rules}
    \minimize_{x \in \mathcal{X}(\mathcal{Z}_r)}
        \sum_{i \in \mathcal{Z}^d}\sum_{j \in \mathcal{Z}^d}
        \Bigl(\theta_{ij}^d - \sum_{k \in \mathcal{K}} P_k^d M_{ij}^k \Bigr) x_{ij},
\end{equation}
where $\mathcal{X}(\mathcal{Z}_r)$ is the feasible set defined by
\eqref{eqn:FOP-Amazon-base-out}--\eqref{eqn:FOP-Amazon-base-restrict}.

\subsubsection{Inverse optimization setup}

We estimate $(\theta^d,P^d)$ from observed zone tours using the suboptimality-based IO
framework of Section~\ref{sec:inverse_milp}. For each route $r \in \mathcal{R}^d$, the
observed zone tour $\tilde{x}$ is treated as approximately optimal for
\eqref{eqn:FOP-Amazon-rules}. As in the production-planning and shift assignment case studies, we enforce the IO optimality
conditions via a cutting-plane procedure: at each iteration, we (i) solve the forward problem
to generate a set of violating alternative tours and (ii) add the corresponding cut to the master
problem.

Operationally, we learn the decision-rule rewards first (fixing $\theta^d$ to its haversine
initialization) and then refine pairwise costs $\theta^d$ on top of the selected rule set.
This two-stage approach isolates broad structural effects (clustering) before fitting
residual pairwise preferences, which is especially important under data sparsity at the
zone-pair level.

\subsubsection{Selection of clustering decision rules}

We next quantify which clustering rules are supported by the data. We evaluate candidate
rules by learning only the corresponding rewards $P^d$ (with $\theta^d$ fixed to haversine
distance) and comparing the summed objective values across depots. At the first level (L1),
we consider four candidate cluster definitions that share three out of the four code
components. The best-performing rule clusters zones by $(W,x,Z)$ (denoted $W\!-\!x.Z$),
achieving the lowest summed objective:
\begin{table}[H]
\centering
\begin{tabular}{ccccc}
\hline
\textbf{L1 cluster rule} & $W\!-\!x.y$ & $W\!-\!x.Z$ & $W.yZ$ & $x.yZ$\\
\hline
\textbf{Objective (summed)} & 160.188 & \textbf{147.144} & 160.855 & 161.093 \\
\hline
\end{tabular}
\caption{Objective values for candidate first-level clustering rules (lower is better).}
\label{tab:L1_cluster}
\end{table}

Conditioned on the selected L1 rule, we evaluate L2 cluster rules by grouping L1 clusters
sharing two of the three components. The best-performing L2 rule clusters by $(W,x)$
(denoted $W\!-\!x$), yielding the lowest summed objective and being best for 15 out of 17
depots:
\begin{table}[H]
\centering
\begin{tabular}{cccc}
\hline
\textbf{L2 cluster rule} & $W\!-\!x$ & $W\!-\!Z$ & $x.Z$\\
\hline
\textbf{Objective (summed)} & \textbf{144.386} & 146.867 & 147.121 \\
\hline
\end{tabular}
\caption{Objective values for candidate second-level clustering rules (lower is better).}
\label{tab:L2_cluster}
\end{table}

We also tested third-level (L3) rules (single-component clustering). These yielded negligible
objective improvement and inconsistent depot-level effects, suggesting weak predictive value
and potential overfitting:
\begin{table}[H]
\centering
\begin{tabular}{ccc}
\hline
\textbf{L3 cluster rule} & $W$ & $x$\\
\hline
\textbf{Objective (summed)} & 144.321 & 144.376 \\
\hline
\end{tabular}
\caption{Objective values for candidate third-level clustering rules (lower is better).}
\label{tab:L3_cluster}
\end{table}

These results support using two clustering levels (L1 and L2) as decision rules and discarding
L3 rules.

\subsubsection{Prediction for test routes and Amazon Score}

For each test route, the input is the set of stops and their associated zones, yielding $\mathcal{Z}_r$.
We first solve the zone-level forward problem \eqref{eqn:FOP-Amazon-rules} with learned
parameters to obtain a predicted zone order. We then enforce this zone order at the stop
level by modifying stop-to-stop travel times:
\[
\tilde{c}_{uv} :=
\begin{cases}
c_{uv}, & \text{if } \mathrm{zone}(u) = \mathrm{zone}(v), \\
c_{uv} + R, & \text{if } \mathrm{zone}(v) \text{ is the next zone after }
                \mathrm{zone}(u) \text{ in the predicted zone tour,} \\
c_{uv} + 2R, & \text{otherwise},
\end{cases}
\]
where $c_{uv}$ is the travel time provided in the dataset and $R$ is a large penalty.
Finally, we solve a stop-level TSP using $\tilde{c}_{uv}$ to obtain the predicted stop
sequence.

Submissions are evaluated using the Amazon Score \citep{noauthor_rc-cliscoring_nodate},
a normalized similarity measure combining sequence deviation ($SD$) and edit distance with
real penalty (ERP):
\[
    \mathrm{Amazon Score}(A,B)
        = \frac{SD(A,B) \cdot ERP_{\mathrm{norm}}(A,B)}{ERP_e(A,B)},
\]
where $A$ is the observed route and $B$ is the predicted route. Lower values indicate
more accurate predictions, with values near zero corresponding to near-exact matches.

\subsubsection{Results}

We report three quantitative findings: (i) the marginal value of decision-rule depth;
(ii) the additional value of pairwise cost learning on top of decision rules; and
(iii) a structural diagnostic based on cluster crossings that explains why decision
rules matter.

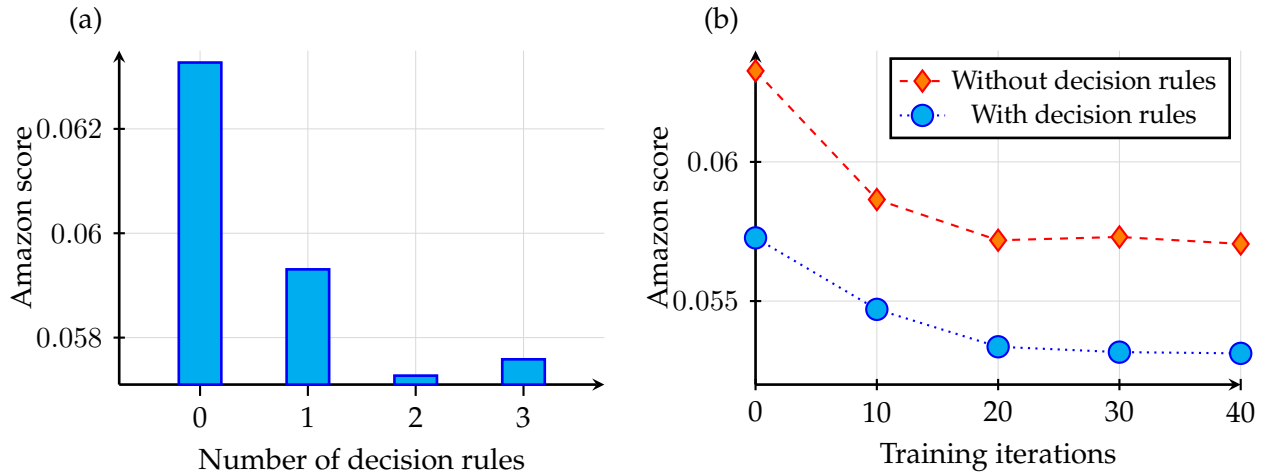
\begin{figure}[!htbp]
    \centering
    \begin{tikzpicture}
    \pgfplotsset{
        myaxis/.style={
            width=8cm,
            height=6cm,
            grid=both,
            grid style={gray!20},
            axis lines=left,
            xlabel style={font=\normalsize},
            ylabel style={font=\normalsize, yshift=-6pt},
            tick style={black, thick},
            tick label style={font=\normalsize},
            legend style={
                draw=black,
                fill=white,
                font=\small
            },
            line width=1pt,
            scaled y ticks=false,
            yticklabel style={
                /pgf/number format/fixed,
                /pgf/number format/precision=4
            }
        }
    }

    \begin{groupplot}[
        group style={
            group name=myplots,
            group size=2 by 1,
            horizontal sep=2cm
        },
        myaxis
    ]

    \nextgroupplot[
        ybar,
        bar width=16pt,
        ymin=0.0571, ymax=0.0635,
        xlabel={Number of decision rules},
        ylabel={Amazon score},
        symbolic x coords={0,1,2,3},
        xtick=data,
        enlarge x limits=0.25,
        major grid style={line width=0.2pt, draw=gray!30},
        minor grid style={line width=0.1pt, draw=gray!10},
        clip=false
    ]
        \refstepcounter{subfigure}
        \node[anchor=south east, font=\normalsize]
            at ([xshift=-2pt,yshift=2pt]current axis.north west)
            {(\alph{subfigure})};

        \addplot[
            fill=cyan,
            draw=blue
        ] coordinates {
            (0,0.0632709)
            (1,0.0593079)
            (2,0.0572715)
            (3,0.0575841)
        };

    \nextgroupplot[
        xmin=0, xmax=4,
        ymin=0.052, ymax=0.064,
        xlabel={Training iterations},
        ylabel={Amazon score},
        xtick={0,1,2,3,4},
        xticklabels={0,10,20,30,40},
        legend style={
            at={(0.98,0.98)},
            anchor=north east,
            draw=black,
            fill=white,
            font=\small
        },
        major grid style={line width=0.2pt, draw=gray!30},
        minor grid style={line width=0.1pt, draw=gray!10},
        clip=false
    ]
        \refstepcounter{subfigure}
        \node[anchor=south east, font=\normalsize]
            at ([xshift=-2pt,yshift=2pt]current axis.north west)
            {(\alph{subfigure})};

        \addplot[
            red,
            dashed,
            thick,
            mark=diamond*,
            mark size=4pt,
            mark options={solid, fill=orange}
        ] coordinates {
            (0,0.0632709)
            (1,0.0586486)
            (2,0.0571852)
            (3,0.0573039)
            (4,0.0570544)
        };
        \addlegendentry{Without decision rules}

        \addplot[
            blue,
            dotted,
            thick,
            mark=*,
            mark size=4pt,
            mark options={solid, fill=cyan}
        ] coordinates {
            (0,0.0572715)
            (1,0.0547026)
            (2,0.0533585)
            (3,0.0531660)
            (4,0.0531218)
        };
        \addlegendentry{With decision rules}

    \end{groupplot}
    \end{tikzpicture}
    \caption{Amazon score as a function of the number of decision rules and training iterations.}
    \label{fig:amazon_combined}
\end{figure}

\paragraph{Effect of decision-rule depth.}
Figure~\ref{fig:amazon_combined} (a) reports the Amazon Score when we include an
increasing number of decision rules in the model (0 = none, 1 = L1 only, 2 = L1+L2,
3 = L1+L2+L3). Starting from a purely distance-based model, adding L1 reduces the score
from $0.0632$ to $0.0593$. Adding L2 yields a further improvement to $0.0573$. In contrast,
adding the L3 rule degrades performance slightly to $0.0576$, consistent with the weak and
likely spurious objective improvements observed in Table~\ref{tab:L3_cluster}.

\paragraph{Decision rules versus pairwise cost learning.}
We next quantify how much improvement is achieved by learning pairwise costs $\theta^d$
after learning decision-rule rewards. Figure~\ref{fig:amazon_combined} (b) reports Amazon Scores
as a function of training iterations (cuts added), comparing models with and without
decision rules.

With decision rules (L1+L2) and \emph{no} pairwise learning, the score is $0.0573$.
As we run IO and refine $\theta^d$, performance improves to $0.0547$ (10 iterations),
$0.0535$ (20 iterations), and stabilizes at $0.0531$ by 30--40 iterations. This corresponds
to a reduction from $0.0573$ to $0.0533$, i.e., an additional $\approx 7.0\%$ improvement
relative to the decision-rule-only model.

Without decision rules, the model starts at $0.0633$ and improves to $0.0586$ (10),
$0.0572$ (20), and stabilizes near $0.0571$ by 40 iterations. Critically,  after
pairwise learning, the no-rules model ($0.0571$) performs similar to the decision-rule-only model
($0.0573$ at iteration 0), with no further pair-wise costs learned, but the fully trained with-rules model achieves $0.0533$. Comparing the
best achieved scores, incorporating decision rules yields a reduction from $0.0571$ to
$0.0533$, i.e., an improvement of approximately $6.7\%$ relative to the best no-rules model,
and approximately $15.8\%$ relative to the distance-based baseline ($0.0633$ to $0.0533$).

\paragraph{Cluster crossings as a structural diagnostic.}
To explain \emph{why} decision rules improve predictive accuracy, we compare the number of
cluster crossings in predicted routes to those in observed driver routes. For each route,
we compute the number of transitions that move from one cluster to another, at L1 and L2
levels, and report these as a percentage of the corresponding count in the observed route.

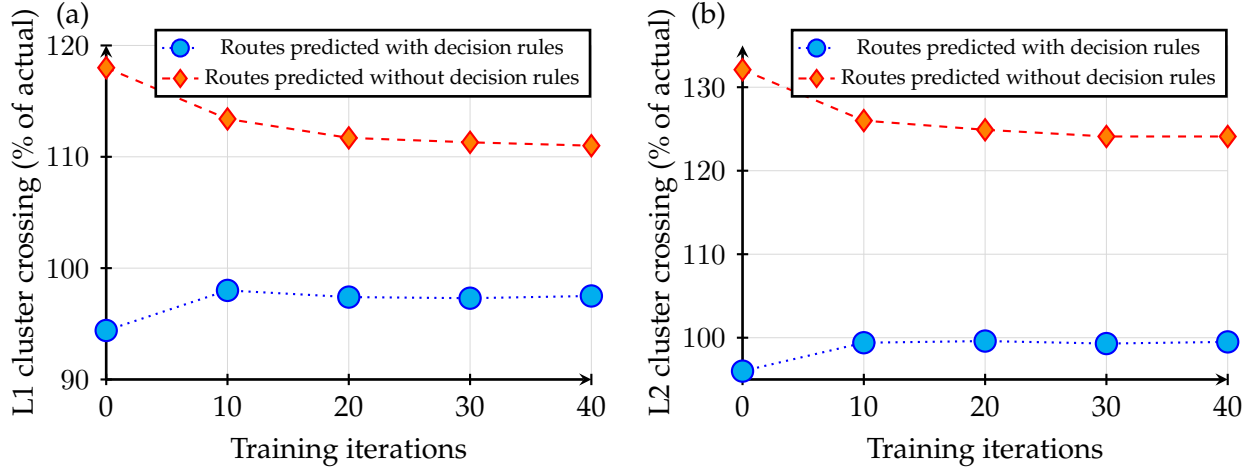
\begin{figure}[!htbp]
    \centering
    \begin{tikzpicture}
    \pgfplotsset{
        myaxis/.style={
            width=8cm,
            height=6cm,
            grid=both,
            grid style={gray!20},
            axis lines=left,
            xlabel style={font=\normalsize},
            ylabel style={font=\normalsize, yshift=-6pt},
            tick style={black, thick},
            tick label style={font=\normalsize},
            legend style={
                draw=black,
                fill=white,
                font=\footnotesize,
                row sep=1pt,
                inner xsep=4pt,
                inner ysep=2pt
            },
            legend image post style={scale=0.8},
            line width=1pt,
            scaled y ticks=false,
            yticklabel style={
                /pgf/number format/fixed,
                /pgf/number format/precision=4
            }
        }
    }

    \begin{groupplot}[
        group style={
            group name=myplots,
            group size=2 by 1,
            horizontal sep=2cm
        },
        myaxis
    ]

    \nextgroupplot[
        xmin=0, xmax=4,
        ymin=90, ymax=120,
        xlabel={Training iterations},
        ylabel={L1 cluster crossing (\% of actual)},
        xtick={0,1,2,3,4},
        xticklabels={0,10,20,30,40},
        legend style={
            at={(1,1.05)},
            anchor=north east,
            draw=black,
            fill=white,
            font=\scriptsize,
            row sep=0pt,
            inner xsep=2pt,
            inner ysep=1pt
        },
        major grid style={line width=0.2pt,draw=gray!30},
        minor grid style={line width=0.1pt,draw=gray!10},
        clip=false
    ]
        \refstepcounter{subfigure}
        \node[anchor=south east, font=\normalsize]
            at ([xshift=-2pt,yshift=2pt]current axis.north west)
            {(\alph{subfigure})};

        \addplot[
            blue,
            dotted,
            thick,
            mark=*,
            mark size=4pt,
            mark options={solid, fill=cyan}
        ] coordinates {
            (0,0.944*100)
            (1,0.980*100)
            (2,0.974*100)
            (3,0.973*100)
            (4,0.975*100)
        };
        \addlegendentry{Routes predicted with decision rules}

        \addplot[
            red,
            dashed,
            thick,
            mark=diamond*,
            mark size=4pt,
            mark options={solid, fill=orange}
        ] coordinates {
            (0,1.180*100)
            (1,1.134*100)
            (2,1.117*100)
            (3,1.113*100)
            (4,1.110*100)
        };
        \addlegendentry{Routes predicted without decision rules}

    \nextgroupplot[
        xmin=0, xmax=4,
        ymin=95, ymax=135,
        xlabel={Training iterations},
        ylabel={L2 cluster crossing (\% of actual)},
        xtick={0,1,2,3,4},
        xticklabels={0,10,20,30,40},
        legend style={
            at={(1.00,1.05)},
            anchor=north east,
            draw=black,
            fill=white,
            font=\scriptsize,
            row sep=0pt,
            inner xsep=2pt,
            inner ysep=1pt        
        },
        major grid style={line width=0.2pt,draw=gray!30},
        minor grid style={line width=0.1pt,draw=gray!10},
        clip=false
    ]
        \refstepcounter{subfigure}
        \node[anchor=south east, font=\normalsize]
            at ([xshift=-2pt,yshift=2pt]current axis.north west)
            {(\alph{subfigure})};

        \addplot[
            blue,
            dotted,
            thick,
            mark=*,
            mark size=4pt,
            mark options={solid, fill=cyan}
        ] coordinates {
            (0,0.960*100)
            (1,0.994*100)
            (2,0.996*100)
            (3,0.993*100)
            (4,0.995*100)
        };
        \addlegendentry{Routes predicted with decision rules}

        \addplot[
            red,
            dashed,
            thick,
            mark=diamond*,
            mark size=4pt,
            mark options={solid, fill=orange}
        ] coordinates {
            (0,1.321*100)
            (1,1.260*100)
            (2,1.249*100)
            (3,1.241*100)
            (4,1.241*100)
        };
        \addlegendentry{Routes predicted without decision rules}

    \end{groupplot}
    \end{tikzpicture}
    \caption{L1 and L2 cluster crossing as percentages of the actual route.}
    \label{fig:cluster_crossing_combined}
\end{figure}

Figure~\ref{fig:cluster_crossing_combined} (a) shows that predicted routes with decision rules exhibit
L1 crossing rates close to the observed routes (between $94.4\%$ and $98.0\%$ across
iterations), whereas predicted routes without decision rules cross L1 clusters far more
often (between $111.0\%$ and $118.0\%$). The same effect is amplified at the L2 level:
Figure~\ref{fig:cluster_crossing_combined} (b) shows with-rules predictions remain close to observed
routes (between $96.0\%$ and $99.6\%$), while no-rules predictions exhibit substantially
more L2 crossings (between $124.1\%$ and $132.1\%$). These results show that decision rules
correctly encode the dominant structural regularity in driver behavior, serving zones in
blocks, whereas pairwise cost learning alone cannot reliably reproduce this behavior.

\subsubsection{Interpretability via route-level analysis}

While aggregate metrics such as the Amazon score quantify predictive performance, they do not directly reveal how decision rules influence the structure of the predicted routes. To complement the quantitative results, we examine representative test routes predicted by our model, with and without decision rules. These examples illustrate how the learned rules shape the routing behavior and help diagnose the structural limitations of cost-only inverse optimization.

\begin{figure}[h!]
    \centering
    \subfloat[][\centering Observed route]{\includegraphics[width=0.45\textwidth]{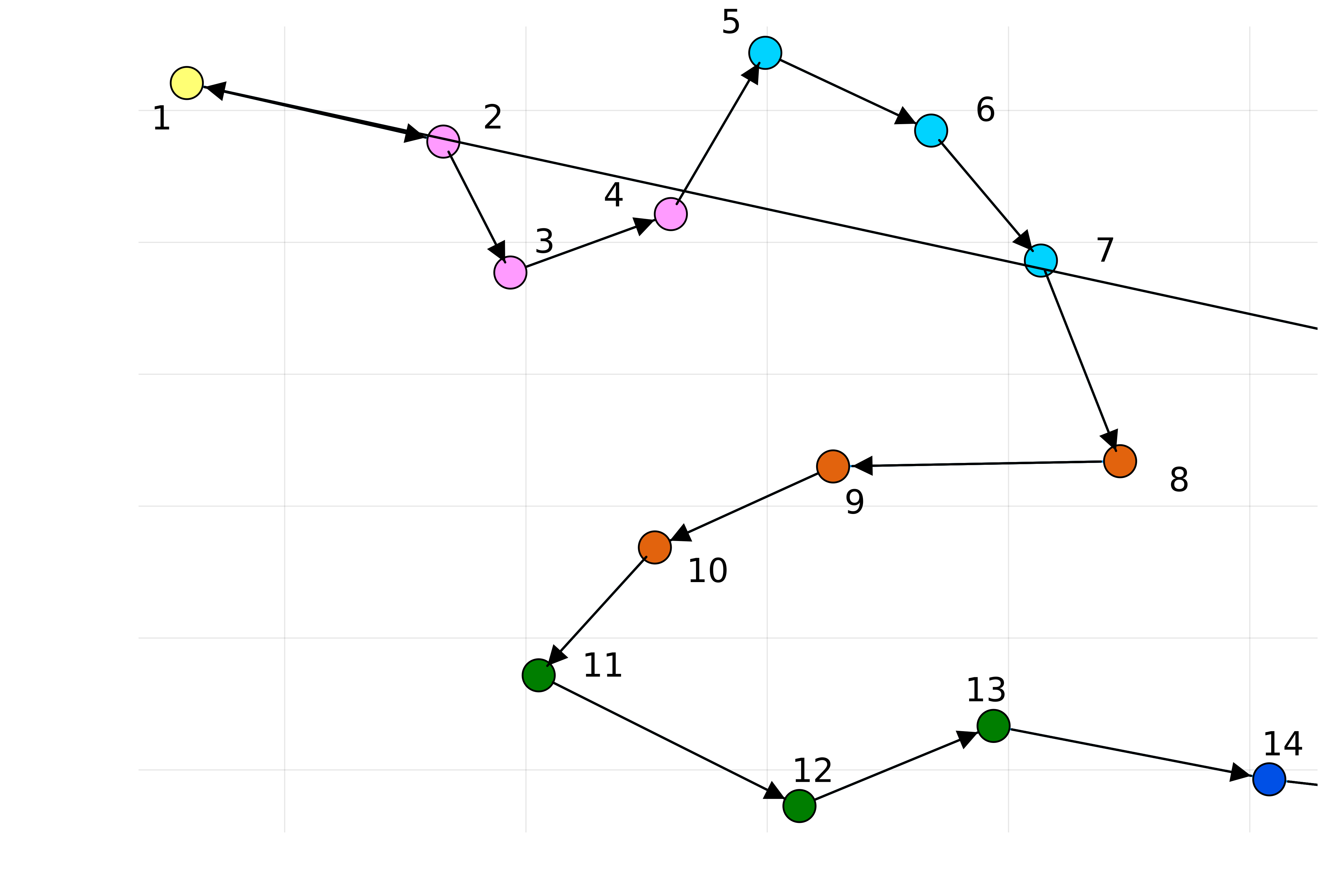}}
    \subfloat[][\centering Route predicted using TSP with distances]{\includegraphics[width=0.45\textwidth]{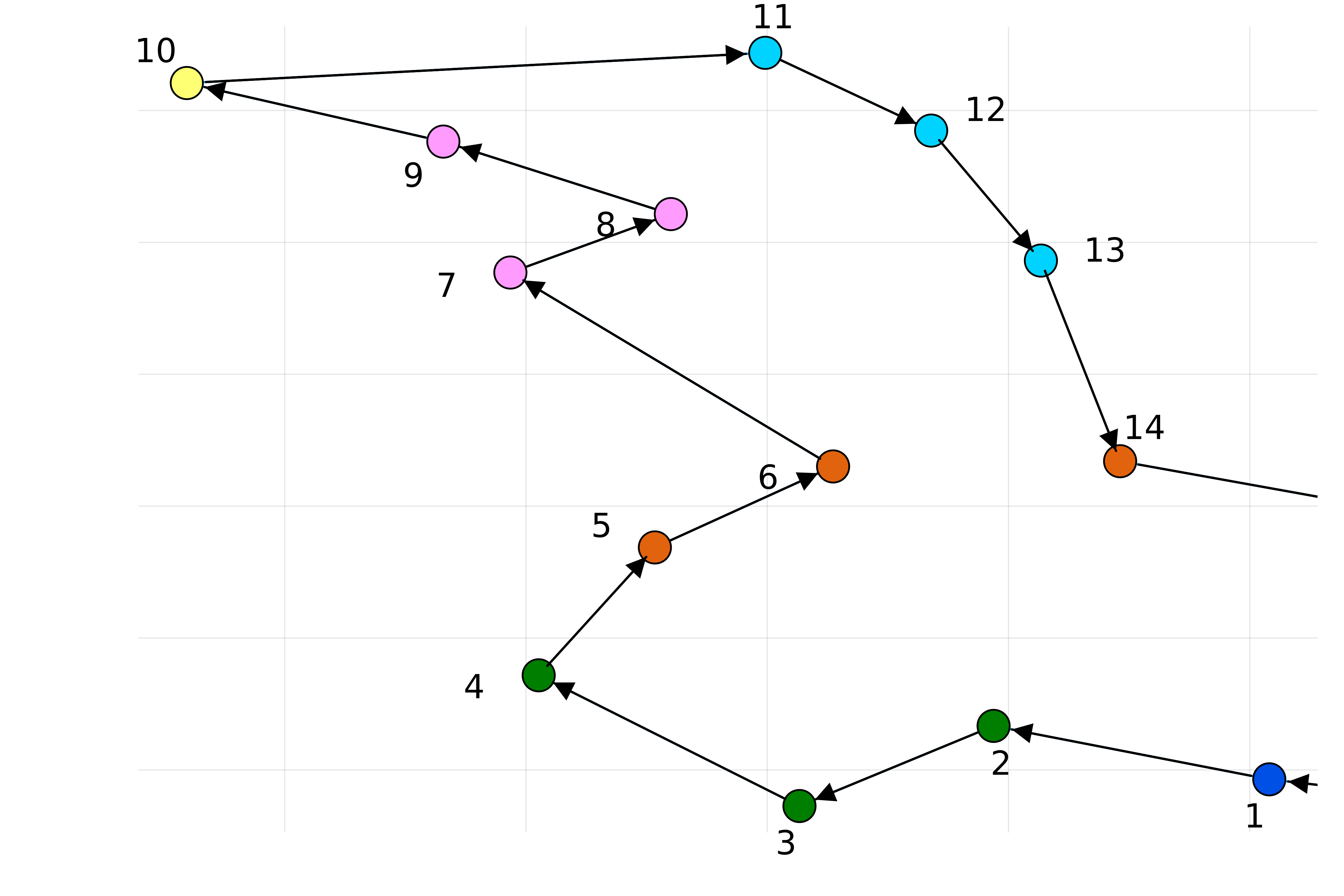}}
    \quad
    \subfloat[][\centering Route predicted using IO w/o decision rules]{\includegraphics[width=0.45\textwidth]{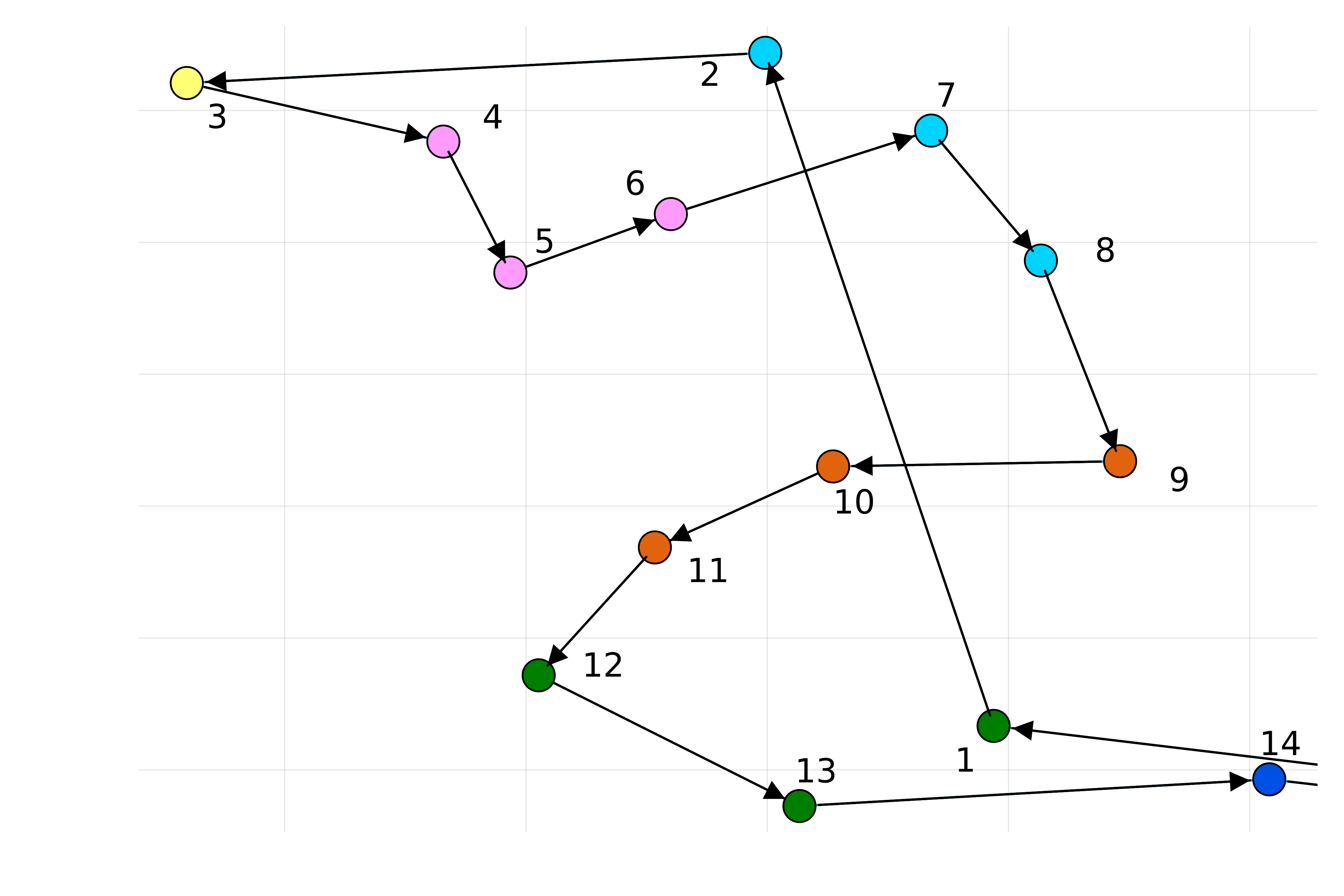}}
    \subfloat[][\centering Route predicted using IO with decision rules]{\includegraphics[width=0.45\textwidth]{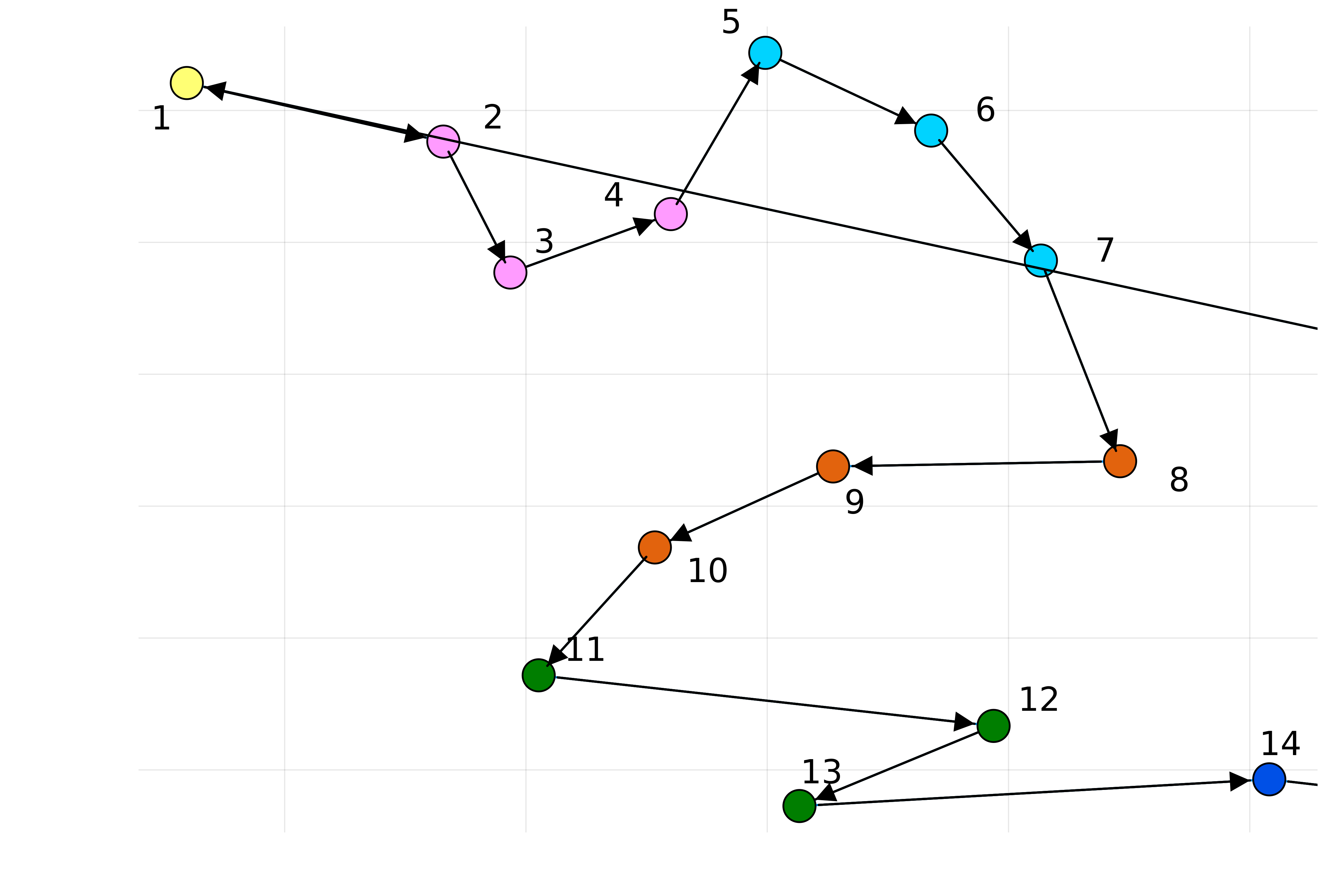}}
    \caption{Actual and predicted routes: Example 1}
    \label{fig:L1_paths_1}
\end{figure}

\begin{figure}[h!]
    \centering
    \subfloat[][\centering Observed route]{\includegraphics[width=0.45\textwidth]{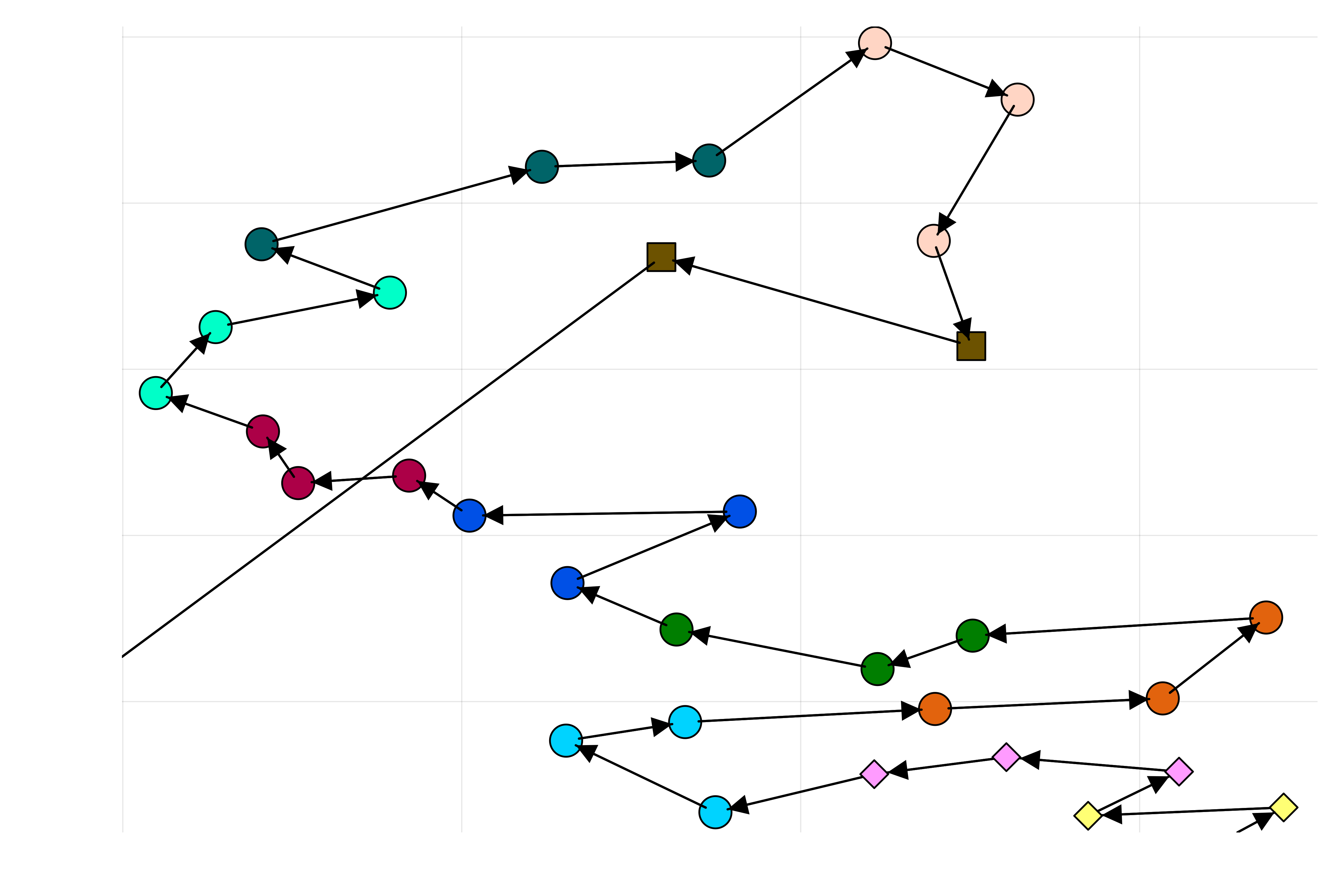}}
    \subfloat[][\centering Predicted route using IO w/o decision rules]{\includegraphics[width=0.45\textwidth]{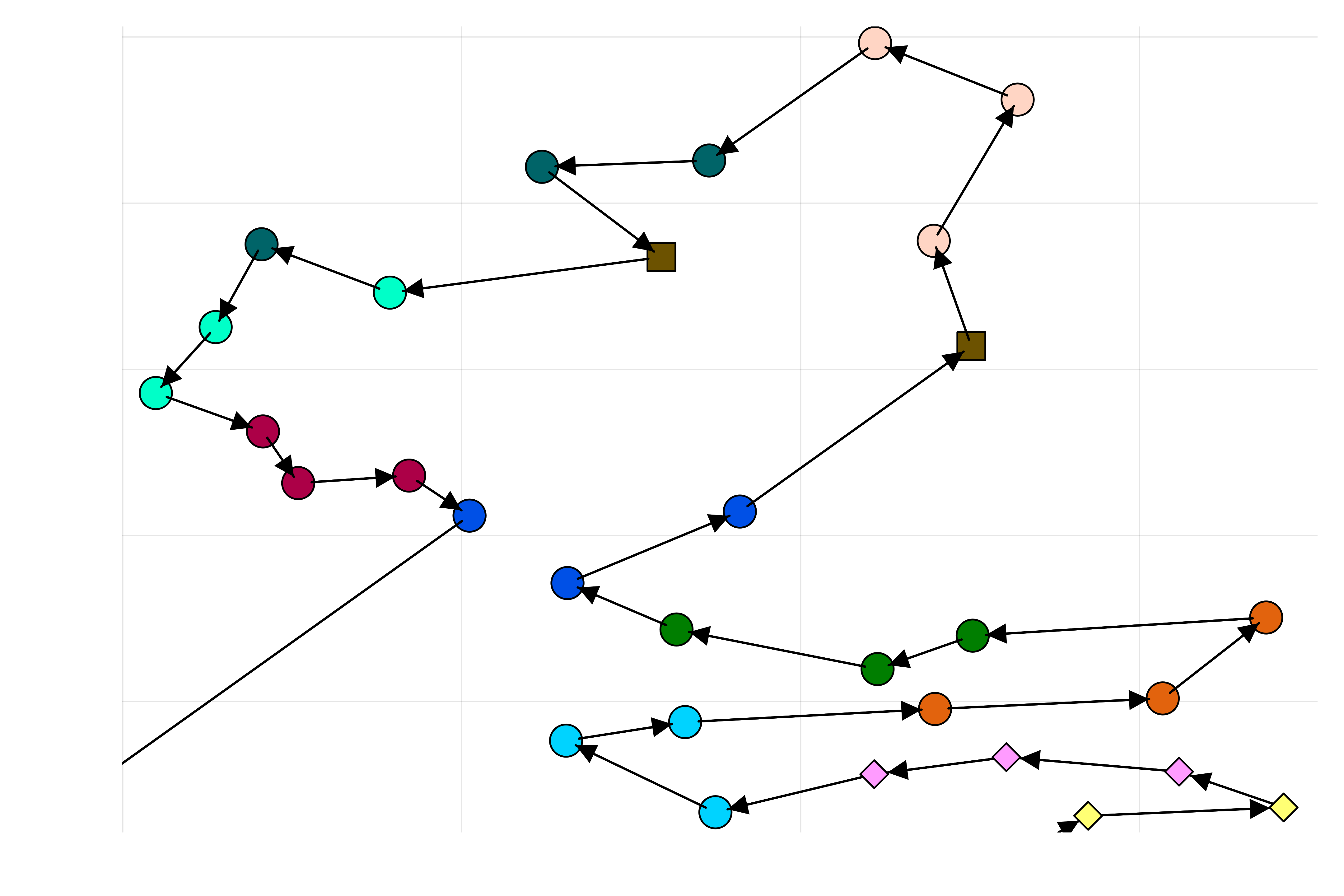}}
    \quad
    \subfloat[][\centering Predicted route using IO with L1 clusters]{\includegraphics[width=0.45\textwidth]{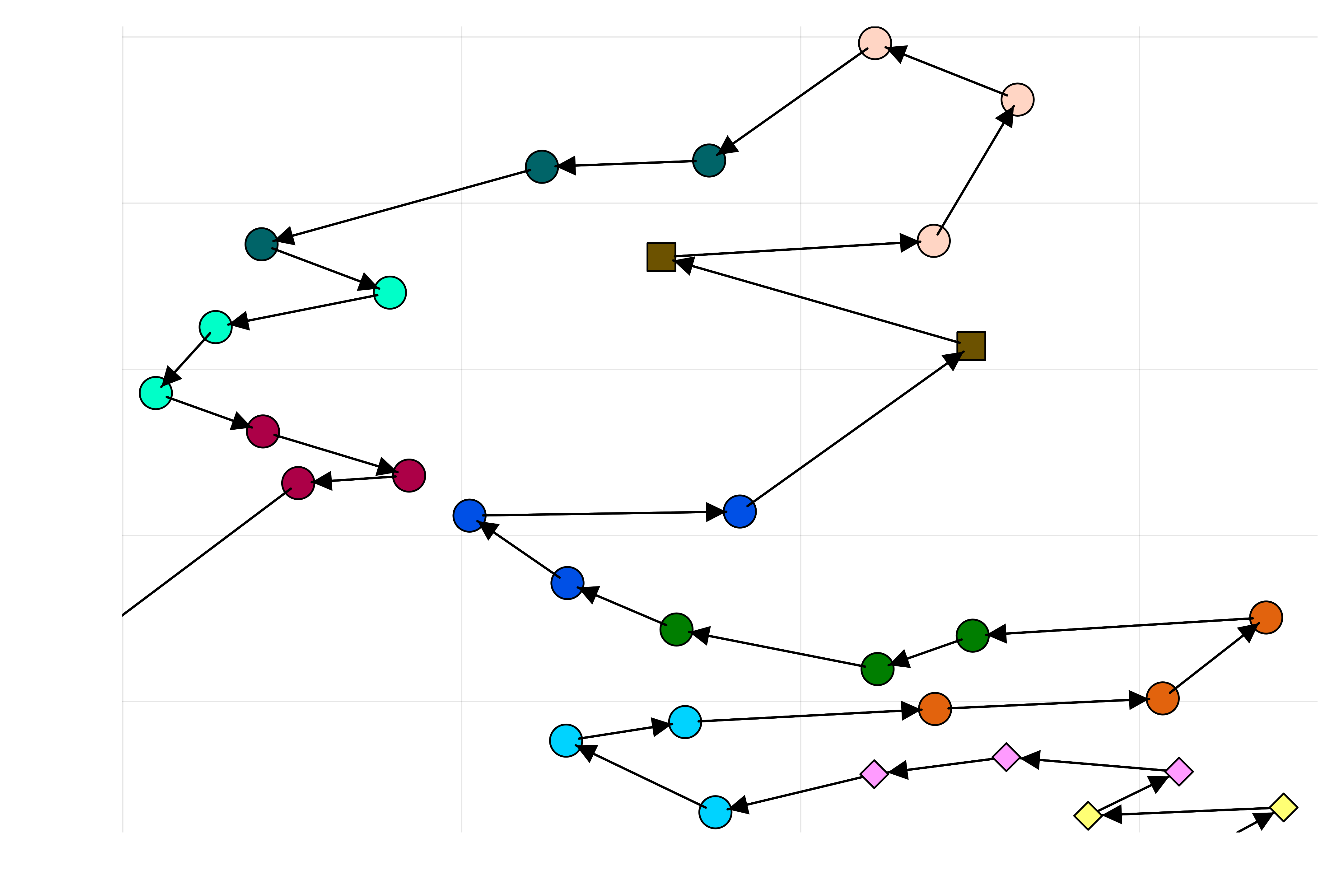}}
    \subfloat[][\centering Predicted route using IO with L2 clusters]{\includegraphics[width=0.45\textwidth]{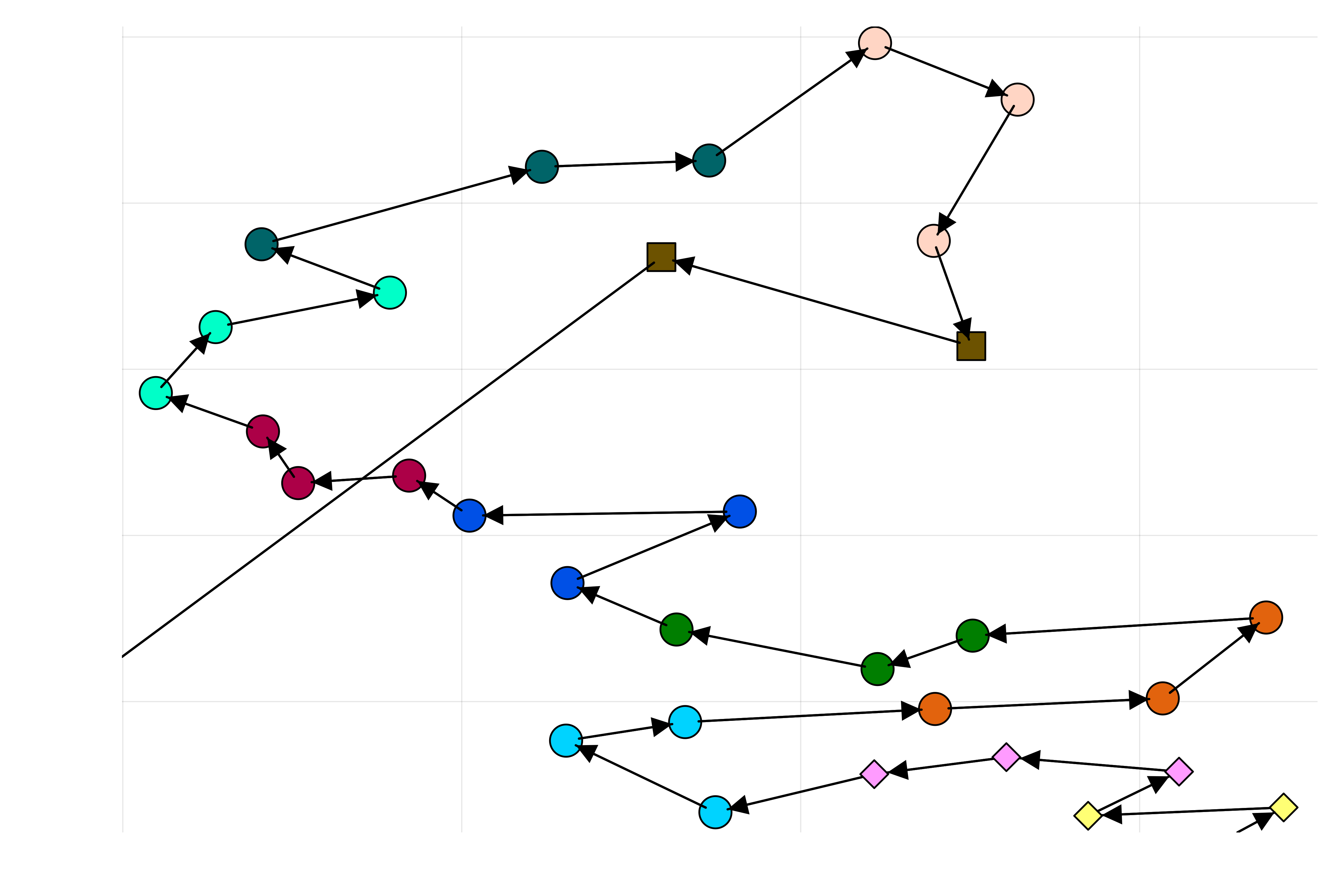}}        
    \caption{Actual and predicted routes: Example 2}
    \label{fig:L1_paths_2}
\end{figure}

All routes shown in Figures~\ref{fig:L1_paths_1} and~\ref{fig:L1_paths_2} are generated by our implementation using the public Amazon Challenge dataset. We study the routes predicted with and without the use of decision rules, to see the impact which learning these decision rule has on the overall accuracy of our route predictions, and discuss how learning the clustering of zones helps improve the final predictions.
% In these figures, zone-level travel costs are based solely on Haversine distances, and the only behavioral signal incorporated is the presence or absence of learned decision rules (i.e., the pairwise cost parameters $\theta_{ij}$ are not yet trained). This isolates the effect of decision rules from fine-grained cost learning.

We first consider the route shown in Figure~\ref{fig:L1_paths_1}, which consists of six first-level (L1) clusters, indicated by the color-coded zones. The observed driver route (Figure~\ref{fig:L1_paths_1}a) exhibits strong contiguity within each cluster, with transitions between clusters occurring only after all zones within a cluster have been served.
However, when the route is predicted solely on the basis of the distance between the zones (Figure~\ref{fig:L1_paths_1}b), or when IO is used to learn the travel costs between the zones, but with the decision rules omitted (Figure~\ref{fig:L1_paths_1}c), the predicted route interleaves zones from different clusters, leading to substantially more cross-cluster transitions.
Finally, when L1 clustering rules are incorporated into the forward model (Figure~\ref{fig:L1_paths_1}d), the predicted route closely matches the observed route: the zones are visited in the same order as in the observed route, and the clusters are visited contiguously.
This observation illustrates that learning the clustering patterns of observed routes allows us to incorporate higher-level structural patterns into our prediction framework, improving the overall accuracy of our predictions.
% In contrast, when decision rules are omitted (Figure~\ref{fig:L1_paths_1}c), the predicted route interleaves zones from different clusters, leading to substantially more cross-cluster transitions. 
% This behavior reflects the well-known failure mode of cost-only TSP-based models, which optimize local travel costs at the expense of higher-level structural coherence.

To further examine the interaction between multiple decision rules, Figure~\ref{fig:L1_paths_2} presents a route that contains both L1 and second-level (L2) clusters. The L1 clusters are represented by the color of the zones, while the L2 clusters are represented by the shape of the zone marker.
The observed route (Figure~\ref{fig:L1_paths_2}a) respects both clustering levels, serving all zones within the L1 and L2 cluster contiguously. 
As seen in the previous example, using IO without decision rules to predict the route leads to violation of the clustering patterns followed by the drivers (Figure~\ref{fig:L1_paths_2}b).
% When both L1 and L2 rules are enforced (Figure~\ref{fig:L1_paths_2}b), the predicted route preserves this hierarchical structure. 
If only L1 rules are enforced (Figure~\ref{fig:L1_paths_2}c), the model the L1 clusters but partially violates the L2 structure, splitting the second L2 cluster (consisting of circle zone markers) into two segments. 
Finally, when both L1 and L2 rules are enforced (Figure~\ref{fig:L1_paths_2}d), the predicted route preserves the complete hierarchical structure of L1 and L2 clusters seen in the observed route, giving us an exact prediction of the driver's route.

% These examples highlight two key insights. First, decision rules act as soft structural preferences that guide the optimizer toward behaviorally plausible routes. Second, omitting relevant rules leads not only to local violations of the corresponding structure, but can also destabilize other aspects of the route due to compensating cost trade-offs. This mirrors the production planning case study, where missing complexity-control rules caused inverse optimization to recover degenerate objective weights and yield poor predictive performance. Together with the quantitative results, the route-level analysis demonstrates that explicitly modeling decision rules is essential for capturing the hierarchical structure present in real-world routing decisions.

These examples highlight two key insights. First, decision rules act as soft structural preferences that guide the optimizer toward behaviorally plausible routes.
Secondly, it is important to include a sufficiently large set of decision rules, to be able to recover the overall structural patterns seen in the data.
Together with the quantitative results, the route-level analysis demonstrates that explicitly modeling decision rules is essential for capturing the hierarchical structure present in real-world routing decisions.

\paragraph{Summary} The Amazon case study reinforces the main insight from the production-planning example:
when decision rules materially shape observed decisions, explicitly embedding these rules
in the forward model is critical for accurate and interpretable inverse optimization.
Quantitatively, learning two clustering rules (L1+L2) reduces the Amazon Score from $0.0632$
(distance-based baseline) to $0.0573$, and refining pairwise costs via IO further reduces
the score to $0.0533$ after 30--40 iterations. In contrast, learning pairwise costs alone
(without decision rules) saturates near $0.0570$ after 40 iterations. The cluster-crossing
diagnostic provides a behavioral explanation for this gap: decision rules align predicted
routes with the block structure present in the observed driver routes, whereas pairwise
learning alone produces excessive cross-cluster transitions.

\section{Conclusions}
\label{sec:conclusions}

In this work, we have presented a generalized data-driven approach to understanding and modeling the decisions made by human experts, using the framework of IO. By incorporating decision rules within the framework, we are able to learn and account for the additional restrictions or constraints that the decision maker may be considering, in addition to learning the objective function, as traditionally done in the IO literature.
We have showcased the capability of this approach using multiple case studies, highlighting the ability of the proposed framework to learn additional constraints imposed by the decision maker and allow for a more realistic model of the decision maker's behavior.

\setlength{\bibsep}{0pt}
\bibliography{references}

\end{document}